\def\version{1.17}
\def\journal{{\small after submission to ART}}
\def\titlep{Pure states on Cuntz algebras arising 
from geometric progressions}
\documentclass[12pt]{article}
\usepackage{graphicx,ifthen,amssymb,epic,eepic,color}

\newcommand{\qed}{\hbox{\rule[-2pt]{3pt}{6pt}}}
\newcommand{\qedh}{\hfill\qed \\}

\font\germ=eufm10 at12pt

\def\goth#1{\hbox{\germ#1}}

\newcommand{\vv}{\vspace{.3in}}

\newtheorem{Thm}{Theorem}[section]

\newtheorem{rem}[Thm]{Remark}

\newtheorem{ex}[Thm]{Example}
\newtheorem{defi}[Thm]{Definition}
\newtheorem{lem}[Thm]{Lemma}

\newtheorem{prop}[Thm]{Proposition}

\newtheorem{cor}[Thm]{Corollary}

\newcommand{\kn}{\Large\bf
$K\hspace{-.4cm} N$
\Large\bf\vv }
\def\cal#1{\mathcal #1}
\def\con{{\cal O}_{n}}
\def\coni{{\cal O}_{\infty}}
\def\edot{=1,\ldots,n}
\def\pr{{\it Proof.}\quad}

\def\co#1{{\cal O}_{#1}}
\def\ltn{\ell^{2}}
\def\ltno{\ltn_{1}}
\def\disp#1{{\displaystyle #1}}
\addtocounter{footnote}{1}
\def\sftt#1{
\setcounter{equation}{0}
\addtocounter{footnote}{1}
\section{#1}
}

\def\ssft#1{\subsection{#1}}
\def\sssft#1{\subsubsection{#1}}

\begin{document}
%
%
\def\autherp{Katsunori Kawamura}
\def\emailp{e-mail: 
kawamurakk3@gmail.com
}
\def\addressp{{\small {\it College of Science and Engineering, 
Ritsumeikan University,}}\\
{\small {\it 1-1-1 Noji Higashi, Kusatsu, Shiga 525-8577, Japan}}
}
\newcommand{\mline}{\noindent
\thicklines
\setlength{\unitlength}{.1mm}
\begin{picture}(1000,5)
\put(0,0){\line(1,0){1250}}
\end{picture}
\par
 }
\def\ptimes{\otimes_{\varphi}}
\def\delp{\Delta_{\varphi}}
\def\sem{\textsf{M}}
\def\ba{\mbox{\boldmath$a$}}
\def\bb{\mbox{\boldmath$b$}}
\def\bc{\mbox{\boldmath$c$}}
\def\be{\mbox{\boldmath$e$}}
\def\bp{\mbox{\boldmath$p$}}
\def\bq{\mbox{\boldmath$q$}}
\def\bu{\mbox{\boldmath$u$}}
\def\bv{\mbox{\boldmath$v$}}
\def\bw{\mbox{\boldmath$w$}}
\def\bx{\mbox{\boldmath$x$}}
\def\by{\mbox{\boldmath$y$}}
\def\bz{\mbox{\boldmath$z$}}
\def\bY{\mbox{\boldmath$Y$}}
\def\titlepage{

\noindent
{\bf 
\noindent
\thicklines
\setlength{\unitlength}{.1mm}
\begin{picture}(1000,0)(0,-300)
\put(0,0){\kn \knn\, for \journal\, Ver.\,\version}
\put(0,-50){\today,\quad {\rm file:} {\rm {\small \textsf{tit01.txt,\, J1.tex}}}}
\end{picture}
}
\vspace{-2.5cm}
\quad\\
{\small 
\footnote{
\begin{minipage}[t]{6in}
directory: \textsf{\fileplace}, \\
file: \textsf{\incfile},\, from \startdate
\end{minipage}
}}
\quad\\
\framebox{
\begin{tabular}{ll}
\textsf{Title:} &
\begin{minipage}[t]{4in}
\titlep
\end{minipage}
\\
\textsf{Author:} &\autherp
\end{tabular}
}
{\footnotesize	
\tableofcontents }
}

\def\openone{\mbox{{\rm 1\hspace{-1mm}l}}}
\def\goh{{\goth h}}
\def\gpr{geometric progression\ }
\def\Gpr{Geometric progression\ }
\def\gprs{geometric progression state}
\def\gpre{geometric progression embedding}
%
%
%
\setcounter{section}{0}
\setcounter{footnote}{0}
\setcounter{page}{1}
\pagestyle{plain}

%
%
\title{\titlep}
\author{\autherp\thanks{\emailp}
\\
\addressp}
\date{}
\maketitle
%
%
\begin{abstract}
Let ${\cal O}_n$ denote the Cuntz algebra for $n\geq 2$.
We introduce an embedding $f$ of ${\cal O}_m$
into ${\cal O}_n$ arising from a geometric progression 
of Cuntz generaters of ${\cal O}_n$.
By identifying ${\cal O}_m$ with $f({\cal O}_m)$,
we extend Cuntz states on ${\cal O}_m$
to ${\cal O}_n$.
We show (i) a necessary and sufficient condition of
the uniqueness of the extension,
(ii) the complete classification of all such extensions up to
unitary equivalence of their GNS representations, and 
(iii) the decomposition formula of a mixing state
into a convex hull of pure states.
The complete set of invariants of all GNS representations 
by such pure states is given as a certain set of complex unit vectors.
\end{abstract}


%
%
\sftt{Introduction}
\label{section:first}
The aim of this paper is to classify a certain class of pure states
on Cuntz algebras in succession to the previous work \cite{GP08}.
For a unital C$^*$-algebra $A$
and a unital C$^*$-subalgebra $B$ of $A$,
any state $\omega$ on $B$ has an 
{\it extension} $\tilde{\omega}$ on $A$,
that is,  $\tilde{\omega}$ is a state on $A$ which satisfies
$\tilde{\omega}|_B=\omega$
(\cite{Dixmier}, 2.10.1),
but it is not unique in general.
In this paper,
we completely classify extensions of a certain class of pure states 
on Cuntz algebras up to unitary equivalence of their
Gel'fand-Naimark-Segal (=GNS) representations.
In consequence,
a new class of pure states on Cuntz algebras and 
the complete set of their invariants are given.
In this section, we show our motivation, definitions and main theorems.
Their proofs will be given in $\S$ \ref{section:third}.

%
%
\ssft{Motivation}
\label{subsection:firstone}
%
%
\sssft{Classification problem of pure states on Cuntz algebras}
\label{subsubsection:firstoneone}
A central problem of representation theory of groups is
the understanding irreducible representations \cite{KV,Kobayashi}.
For example, 
it means  
construction of irreducible representations,
finding a complete set of  invariants of representations, and 
understanding these invariants.
Our purpose is to study irreducible representations of C$^*$-algebras
according to such subjects.
By GNS construction, 
the state theory of a C$^*$-algebra $A$ can be interpreted as  
the  (cyclic) representation theory of $A$ almost all.
Hence we mainly consider (pure) states instead of (irreducible) 
representations in this paper.

For Cuntz algebras which are
typical examples of separable infinite simple C$^*$-algebras 
(see $\S$ \ref{subsection:firsttwo}),
representations and states have been studied by many authors
\cite{ACE,BC,BJKW,BJO,BK,DaPi2,DaPi3,DHJ,Evans,Gabriel,Jeong1999,
Jeong2005,GP08,Laca1993,LS,Shin}.
They have various applications,
for example,
endomorphisms of ${\cal B}({\cal H})$ \cite{BJ1997,BJP,FL},
iterated function systems \cite{BJ},
Markov measures \cite{DJ14,DJ},
wavelets \cite{JorgensenW},
continued fractions \cite{CFR02},
construction of $R$-matrices \cite{TS15},
construction of multiplicative isometries \cite{TS08},
invariant measures \cite{PFO01},
and string theory \cite{AK03}.
But their classifications have not been finished yet.

We intend to develop the classification of pure states on Cuntz algebras
by constructing a new class of states.
For this purpose,
we take notice of Cuntz states which 
are completely classified pure states with explicit numerical invariants 
(see $\S$ \ref{subsubsection:firsttwotwo}.) 
About other result of complete classification of states, see \cite{TS11}.

As a method to construct new states,
we consider extensions of  Cuntz states in this paper.
A new idea is as follows:
In the previous work \cite{GP08},
we classified extensions of Cuntz states on $\co{n^m}$
to $\con$ for any integer $m\geq 2$
with respect to a certain embedding of $\co{n^m}$
into $\con$.
In this paper, we generalize a method to extend Cuntz states
{\it but not} replace Cuntz states with general states.
This is a crucial point to make a computable theory.
For a general embedding $f$ of $\co{m}$ into $\con$,
we introduce a new notion ``$f$-sub-Cuntz state"
as an extension of a Cuntz state on $\co{m}$
to $\con$ ($\S$ \ref{subsubsection:firsttwotwo}).
As examples of $f$-sub-Cuntz state,
we choose a certain class of embeddings (= \gpre s)
and classify $f$-sub-Cuntz states associated with them
(= \gprs s) ($\S$ \ref{subsubsection:firsttwothree}).
We will more closely explain this idea 
and its merits in the next subsection.

%
%
\sssft{Extensions of Cuntz states arising from embeddings}
\label{subsubsection:firstonetwo}
For $2\leq n\leq \infty$, let $\con$ denote the Cuntz algebra. 
The outline of our strategy is as follows:
\begin{enumerate}
\item[Step 1]
Fix a unital embedding $f$ of $\co{m}$ into $\con$ and identify
$\co{m}$ with $f(\co{m})$.
\item[Step 2]
Extend Cuntz states on $\co{m}$ to $\con$
with respect to the inclusion $\co{m}\hookrightarrow \con$ in Step 1.
\item[Step 3]
Study such extensions:
\begin{enumerate}
\item
For a given Cuntz state $\omega$, 
is an extension of $\omega$ unique? 
\item
If it is unique, then write down its state values explicitly.
\item
Find the condition of equivalence between two extensions.
Furthermore,
find the complete set of invariants of extensions.
\item 
If a parametrization of such (equivalence classes of) extensions are given,
then investigate properties of the parametrization closely.
\end{enumerate}
\end{enumerate}

\noindent
These will be explicitly explained in $\S$ \ref{subsubsection:firsttwotwo} again.
Merits of this scheme are as follows:
\begin{enumerate}
\item
Cuntz states are well studied and they have a good parametrization.
Hence it is expected that 
their generalizations also have good properties.
For example, sub-Cuntz states are successful generalizations \cite{GP08}.
\item
If one succeeds at the proof of the uniqueness of extension,
then 
the purity of the extended state holds automatically (\cite{Ped}, 4.1.7).
\item
Embeddings are available to study states as new tools.
If one chooses adequate embeddings, then 
it is expected that states arising from them
are computable and one can get the complete classification of them.
Furthermore,  one can easily generalize known theorems
by generalizing related embeddings (see Definition \ref{defi:lexico}
and (\ref{eqn:ftisj})).
\item
It is considered that 
this method is a kind of induced representation theory 
in the broad sense of the term. Well-known Rieffel's induction \cite{Rieffel}
requires a conditional expectation (or its generalization) from an algebra 
to a subalgebra in order to define an induced representation.
However, to find a conditional expectation is not easy 
except a few typical classes.
Even if one does not know a conditional expectation,
extensions of states always exist (Step 2) and 
(special classes of) embeddings can be easily constructed.
Hence we expect that this scheme can play a role of alternatives
of the induced representation theory of C$^*$-algebras.
\end{enumerate}

This paper includes results in \cite{GP06}
by interpreting ``representation" as ``state".

%
%
\ssft{Geometric progression states}
\label{subsection:firsttwo}
In this subsection, 
first, we will introduce a class of states arising 
from a general embedding of Cuntz algebras.
Next, we will define geometric progression states as its special case.
For $2\leq n\leq \infty$, 
let $\con$ denote the {\it Cuntz algebra} \cite{C}, that is, 
$\con$ is a C$^{*}$-algebra 
which is universally generated by a (finite or infinite) sequence
$s_{1},\ldots,s_n$ satisfying
$s_{i}^{*}s_{j}=\delta_{ij}I$ for $i,j\edot$ and
%
%
\begin{equation}
\label{eqn:cuntzdef}
\sum_{i=1}^{n}s_{i}s_{i}^{*}=I\mbox{ when } n<\infty,\quad
\sum_{i=1}^{k}s_{i}s_{i}^{*}\leq I,\quad k= 1,2,\ldots 
\mbox{when }n = \infty
\end{equation}
where $I$ denotes the unit of $\con$.
The Cuntz algebra $\con$ is an infinite dimensional, 
noncommutative C$^{*}$-algebra with unit.
Furthermore, $\con$ is {\it simple}, that is,
there exists no nontrivial closed two-sided ideal of $\con$.

%
%
\sssft{Embeddings of Cuntz algebras}
\label{subsubsection:firsttwoone}
We review basics of general embeddings of Cuntz algebras.
For two unital C$^*$-algebras  $A$ and $B$,
let ${\rm Hom}(A,B)$ denote
the set of all unital $*$-homomorphisms from $A$ to $B$.
If $A$ is simple,
then any $f\in {\rm Hom}(A,B)$ is
injective, that is,
$f$ is an embedding of $A$ into $B$.
In this paper,
we consider only unital embeddings. 
Let $s_1,\ldots,s_n$ denote Cuntz generators of $\con$.
In general,
$f\in {\rm Hom}(\con,A)$ is identified
with Cuntz generators $S_1,\ldots,S_n$ in $A$
as $f(s_i)=S_i$ for $i=1,\ldots,n$.
Hence we can define $f$ by only images $\{f(s_i)\}_{i=1}^n$.
%
%
\begin{rem}
\label{rem:first}
{\rm
(\cite{SE01}, Lemma 2.1)
For $2\leq m,n<\infty$,
${\rm Hom}(\co{m},\con)\ne \emptyset$
if and only if $m=(n-1)k+1$ for some $k\geq 1$.
In this paper,
we always assume the latter condition for $(m,n)$ for a given inclusion
$\co{m}\subset \con$.
}
\end{rem}

%
%
\sssft{$f$-sub-Cuntz states}
\label{subsubsection:firsttwotwo}
In this subsection, we introduce $f$-sub-Cuntz states.
For this purpose, we review Cuntz state on $\con$.
For $2\leq n\leq \infty$,
let $s_1,\ldots,s_n$ denote Cuntz generators of $\con$.
For any complex unit vector $z=(z_1,\ldots,z_n)\in {\Bbb C}^n$,
a state $\omega_z$ on $\con$ which satisfies
%
%
\begin{equation}
\label{eqn:first}
\omega_z(s_j)=\overline{z_j}\quad\mbox{for all }j=1,\ldots,n,
\end{equation}
exists uniquely and is pure,
where $\overline{z_j}$ denotes the complex conjugate of $z_j$.
When $n=\infty$,
replace ${\Bbb C}^n$ by $\ltn:=\{(z_j): \sum_{j\geq 1}|z_j|^2=1\}$.
The state $\omega_z$ is called the {\it Cuntz state} 
by $z$ \cite{BJ1997,BJ,BJP,BK}.
GNS representations 
by $\omega_z$ and $\omega_y$
are unitarily equivalent if and only if $z=y$
(see Appendix B in \cite{GP08}).
Bratteli and Jorgensen \cite{BJ} introduced 
sub-Cuntz states as generalizations of Cuntz states 
(see $\S$ \ref{subsection:secondtwo}).
Furthermore, we generalize sub-Cuntz states as follows. 

Fix $f\in {\rm Hom}({\cal O}_m,{\cal O}_n)$ 
and identify $\co{m}$ with $f(\co{m})\subset \con$
for $2\leq m\leq \infty$.
With respect to this identification,
an extension of a Cuntz state on $\co{m}$
to $\con$ always exists.
We call such a state as an {\it $f$-sub-Cuntz state} on $\con$.
More concretely,
for a unit vector $z=(z_{1},\ldots,z_m)\in {\Bbb C}^m$,
$\omega$ is an $f$-sub-Cuntz state on $\con$ by $z$
if $\omega$ is a state on $\con$ which satisfies the following equations:
%
%
\begin{equation}
\label{eqn:fsub}
\omega(f(t_j))=\overline{z_j}\quad \mbox{for all }j=1,\ldots,m
\end{equation}
where $t_1,\ldots,t_m$ denote Cuntz generators of $\co{m}$.
When $n=\infty$,
replace ${\Bbb C}^n$ by $\ltn$.
The most essential properties of $f$-sub-Cuntz state
are as follows.
%
%
\begin{lem}
\label{lem:homo}
Fix $2\leq m\leq \infty$ and
$f\in {\rm Hom}(\co{m},\con)$.
Let $V:={\Bbb C}^m$ when $m<\infty$ and $V:=\ltn$ when $m=\infty$.
Define $V_1:=\{x\in V:\|x\|=1\}$.
\begin{enumerate}
\item(Existence)
For any $z\in V_1$,
an $f$-sub-Cuntz state on $\con$ by $z$ exists.
If it is unique, then it is pure.
\item(Equivalent conditions)
For $z\in V_1$
and a state $\omega$ on $\con$ with the GNS representation
$({\cal H},\pi,\Omega)$,
the following are equivalent:
\begin{enumerate}
\item
$\omega$ is an $f$-sub-Cuntz state by $z$.
\item
$\disp{\sum_{j}z_j \pi(f(t_j))\Omega=\Omega}$.
\item
$\pi(f(t_j)^*)\Omega=z_j\Omega$ for all $j$.
\end{enumerate}
\end{enumerate}
\end{lem}
%
%
\pr
(i)
The existence follows by definition.
Since any Cuntz state is pure,
there always exists a pure extension (\cite{Dixmier}, 2.10.1).
Especially, if its extension is unique,
then it is pure automatically.

\noindent
(ii)
Let $T_j:=f(t_j)$.

\noindent
(a)$\Rightarrow$(b)
If $\omega$ is an $f$-sub-Cuntz state by $z$,
then $\sum_{j}z_j\omega(T_j)=1$ from (\ref{eqn:fsub}).
This implies $\|\sum_{j}z_j\pi(T_j)\Omega-\Omega\|=0$.
Hence (b) is proved.\\
(b)$\Rightarrow$(c)
By operating $\pi(T_j)^*$ at both sides in (b),
(c) is obtained.\\
(c)$\Rightarrow$(a) 
$\omega(T_j)
=\langle\Omega|\pi(T_j)\Omega\rangle
=\langle \pi(T_j)^*\Omega|\Omega\rangle
=\langle z_j\Omega|\Omega\rangle=\overline{z_j}$.
\qedh

%
%
\begin{rem}
\label{rem:second}
{\rm
\begin{enumerate}
\item
An $f$-sub-Cuntz state by $z$ is not always unique 
(see Theorem \ref{Thm:mainb}(i) or Theorem 
\ref{Thm:maintwo}(ii)).
\item
Any Cuntz state on $\con$
is an $f$-sub-Cuntz state with respect to
$f=id_{\con}$.
\item
In (\ref{eqn:first}),
only special values of an $f$-sub-Cuntz state by $z$ are given,
but not its general values.
Hence the determination of all values of a given $f$-sub-Cuntz state
is one of fundamental problems.
\item
If $n=m$ and $f$ is bijective,
then $f$ is an automorphism of $\con$.
Let $\alpha$ denote the standard $U(n)$-action on $\con$.
If $f=\alpha_g$ for $g\in U(n)$,
then a transformation of any Cuntz state by $f$ is also a Cuntz state.
In general, 
if $f$ is an automorphism of $\con$,
then an $f$-sub-Cuntz state by any $z\in ({\Bbb C}^n)_1$ is unique.
\end{enumerate}
}
\end{rem}

\noindent
We will show other properties of $f$-sub-Cuntz states 
in $\S$ \ref{subsection:secondfour}.

%
%
\sssft{Geometric progression states}
\label{subsubsection:firsttwothree}
In this subsection,
we introduce a special class of $f$-sub-Cuntz states on $\con$.
For $2\leq n<\infty$
and $2\leq m\leq \infty$,
let $s_1,\ldots,s_n$ and $t_1,\ldots,t_m$ denote
Cuntz generators of $\con$ and $\co{m}$, respectively.
%
%
\begin{defi}
\label{defi:lexico}
\begin{enumerate}
\item
When $m=(n-1)k+1$ for $k\geq 1$,
define $f\in {\rm Hom}(\co{m},\con)$ by
%
%
\begin{equation}
\label{eqn:gencc}
\left\{
\begin{array}{rll}
f(t_{(n-1)r+i}):= 
&s_{n}^{r}s_{i}\quad&
\left(
\begin{array}{l}
r=0,1,\ldots,k-1,\\
i=1,\ldots,n-1
\end{array}
\right),\\
&\\
f(t_{m}):= &s_{n}^{k}.
\end{array}
\right.
\end{equation}
%
\item
When $m=\infty$,
define $f\in {\rm Hom}(\coni,\con)$ by
%
%
\begin{equation}
\label{eqn:gencci}
f(t_{(n-1)r+i}):= s_n^{r}s_i\quad (r\geq 0,\,i=1,\ldots,n-1).
\end{equation}
\end{enumerate}

\noindent
Here $s_{n}^{0}$ denotes the unit of $\con$
for convenience.
We call $f$ the {\it geometric progression embedding} of $\co{m}$
into $\con$. 
\end{defi}

\noindent
For example,
when $(n,m)=(2,3)$,
$f$ in (\ref{eqn:gencc}) is given as follows:
%
%
\begin{equation}
\label{eqn:ftone}
f(t_1)=s_1,\quad f(t_2)=s_2s_1,\quad f(t_3)=s_2^2.
\end{equation}
When $n=2$ in (\ref{eqn:gencci}),
$f(t_i)$'s are given as follows:
%
%
\begin{equation}
\label{eqn:geom}
s_1,s_2s_1,s_2^2s_1,s_2^3s_1,\ldots, s_2^rs_1,\ldots.
\end{equation}
This is the origin of ``geometric progression embedding".
Geometric progression embeddings have appeared in 
\cite{C,Davidson,GP06,BFO01,SE01,UFB01,PC01,CFR01}.
We expressly provide
the definition of $f$-sub-Cuntz state
for a given geometric progression embedding $f$ according to 
(\ref{eqn:fsub}) as follows.
\begin{enumerate}
\item
For $f$ in (\ref{eqn:gencc}),
$\omega$ is an $f$-sub-Cuntz state by $z=(z_1,\ldots,z_m)\in ({\Bbb C}^m)_1
:=\{y\in {\Bbb C}^m:\|y\|=1\}$
if and only if $\omega$ is a state on $\con$ which satisfies
%
%
\begin{equation}
\label{eqn:fsubl}
\left\{
\begin{array}{rl}
\omega(s_n^rs_i)=&\overline{z_{(n-1)r+i}}\quad \mbox{for all }
r=0,1,\ldots,k-1 \mbox{ and }
i=1,\ldots,n-1,\\
\\
\omega(s_n^k)=&\overline{z_m}.
\end{array}
\right.
\end{equation}
\item
For $f$ in (\ref{eqn:gencci}),
$\omega$ is an $f$-sub-Cuntz state on $\con$ 
for $z=(z_1,z_2,\ldots)\in\ltno:=\{y\in\ltn:\|y\|=1\}$
if and only if $\omega$ is a state on $\con$ which satisfies
%
%
\begin{equation}
\label{eqn:fsubli}
\omega(s_n^{r}s_i)
=\overline{z_{(n-1)r+i}}\quad \mbox{for all }r\geq 0,\,
i=1,\ldots,n-1.
\end{equation}
\end{enumerate}
We call these {\it 
geometric progression states on $\con$ by $z$ of order $k$ and of order $\infty$}
in (\ref{eqn:fsubl}) and (\ref{eqn:fsubli}), respectively.
For $f$ in (\ref{eqn:gencc}),
when $k=1$, $f=id_{\con}$.
In this case, any $f$-sub-Cuntz state is just a Cuntz state.

About varieties of \gpre s and states arising from them,
see $\S$ \ref{subsection:fourththree}.

%
%
\ssft{Main theorems}
\label{subsection:firstthree}
In this subsection, we show our main theorems.
From Lemma \ref{lem:homo}(i),
remaining problems on $f$-sub-Cuntz states
are the uniqueness,
decomposition formulas of mixture states, 
and their equivalence.
We consider these problems for states 
introduced in $\S$ \ref{subsubsection:firsttwothree}.
%
%
\sssft{Uniqueness, purity and decomposition of mixture}
\label{subsubsection:firstthreeone}
%
%
\begin{Thm}(Uniqueness)
\label{Thm:mainb}
\begin{enumerate}
\item
Let $m=(n-1)k+1$ for $k\geq 2$.
For $z=(z_1,\ldots,z_m)\in ({\Bbb C}^m)_1$,
a geometric progression state $\omega$ on $\con$ by $z$
is unique if and only if $|z_m|<1$.
In this case, 
$\omega$ is pure.
We write $\omega$ as $\omega_z'$ when $|z_m|<1$.
\item
For any $z\in \ltno$,
a geometric progression state on $\con$ by $z$ is unique and pure.
We write such a state as $\omega_z'$.
\end{enumerate}
\end{Thm}
%

\noindent
Remark that if $k=1$ in Theorem \ref{Thm:mainb}(i),
then $m=n$ and $\omega$ is just the Cuntz state
by $z$, which is unique for all $z\in ({\Bbb C}^n)_1$.

%
%
\begin{Thm}
\label{Thm:mixture}
(Decomposition of mixture)
For $m=(n-1)k+1$ with $k\geq 2$ and 
$z=(z_1,\ldots,z_m)\in ({\Bbb C}^m)_1$,
if $|z_m|=1$, then  
a geometric progression state on $\con$ by $z$ is a convex hull
of Cuntz states by $(0,\ldots,0,e^{2\pi j\sqrt{-1}/k}q)\in {\Bbb C}^n$
for $j=1,\ldots,k$
where $q$ is a $k$-th root of $z_m$.
\end{Thm}

\noindent
Since any Cuntz state is a geometric progression state,
the set of geometric progression states is closed with respect to
the pure state decomposition from Theorem \ref{Thm:mixture}.

From Theorem \ref{Thm:mainb}(i) and Theorem \ref{Thm:mixture},
the following holds.
%
%
\begin{cor}
\label{cor:purity}
For $m=(n-1)k+1$ for $k\geq 2$,
a geometric progression state $\omega$ by $z=(z_1,\ldots,z_m)\in
({\Bbb C}^m)_1$
is pure if and only if 
\begin{enumerate}
\item
$|z_m|<1$ or,
\item
$|z_m|=1$  and 
$\omega$ is the Cuntz state
by $(0,\ldots,0,q)\in {\Bbb C}^n$ for some $k$-th root $q$ of $z_m$.
\end{enumerate}
\end{cor}

%
%
\sssft{Equivalence}
\label{subsubsection:firstthreetwo}
For two states $\omega$ and $\omega'$ on $\con$,
we write
$\omega\sim \omega'$ if 
their GNS representations are unitarily equivalent.
%
%
\begin{Thm}
\label{Thm:mainthree}
Let $\omega_z'$ be as in Theorem \ref{Thm:mainb}.
\begin{enumerate}
\item
For $m=(n-1)k+1$ with $k\geq 2$,
define
${\cal W}_m:=\{(w_1,\ldots,w_m)\in ({\Bbb C}^m)_1:|w_m|<1\}$.
For $z,y\in {\cal W}_m$,
$\omega_z'\sim \omega_y'$ if and only if $z=y$.
\item
For $z,y\in \ltno$,
$\omega_z'\sim \omega'_y$ if and only if $z=y$.
\end{enumerate}
\end{Thm}

\noindent
From Theorem \ref{Thm:mainthree},
it is shown that ${\cal W}_m$ ({\it resp.} $\ltno$) 
is the complete set of invariants of unitary equivalence classes 
of pure geometric progression states on $\con$ of order $k$ ({\it resp.}
of order $\infty$).

Next,  we show the equivalence condition
between $\omega_z'$ ($z\in \ltno$) and $\omega_y'$ ($y\in {\cal W}_m$).
%
%
\begin{Thm}
\label{Thm:reductionequiv}
Let $\{e_i\}$ denote the standard basis of $\ltn$ and let $z\in \ltno$.
\begin{enumerate}
\item
For $y=(y_1,\ldots,y_n)\in ({\Bbb C}^n)_1$,
let $\omega_y$ be as in (\ref{eqn:first}).
Then 
$\omega_z'\sim \omega_y$ if and only if 
$|y_n|<1$
and $z=\tilde{y}$
where $\tilde{y}\in \ltno$ is defined as 
%
%
\begin{equation}
\label{eqn:tilde}
\tilde{y}:=\sum_{r\geq 0}
\sum_{i=1}^{n-1}y_{n}^r y_{i}\,e_{(n-1)r+i}.
\end{equation}
In this case,  $\omega_z'=\omega_y$.
\item
Assume $m=(n-1)k+1$ and $k\geq 2$.
For $y=(y_1,\ldots,y_m)\in {\cal W}_{m}$,
$\omega_z'\sim \omega_y'$
if and only if $z=\tilde{y}$
where $\tilde{y}\in \ltno$ is defined as
%
%
\begin{equation}
\label{eqn:tildebb}
\tilde{y}:=\sum_{r\geq 0}
\sum_{i=1}^{m-1}y_{m}^r y_{i}\,e_{(m-1)r+i}.
\end{equation}
In this case, 
$\omega_z'=\omega_y'$.
\end{enumerate}
\end{Thm}

\noindent
From Theorem \ref{Thm:reductionequiv}(i) and (ii),
and Theorem \ref{Thm:mainthree}(ii),
the complete set of invariants of all pure \gprs s on $\con$
is given as follows:
%
%
\begin{equation}
\label{eqn:cset}
\ltno\cup\{(0,\ldots,0,c)\in{\Bbb C}^n:|c|=1\}.
\end{equation}
In other words,
every pure \gprs\ is parametrized by a vector in (\ref{eqn:cset}),
and
for any two distinct vectors in (\ref{eqn:cset}),
associated \gprs s are not equivalent.

Next, we show relations between \gprs s of different finite orders.
%
%
\begin{Thm}
\label{Thm:mequivalence}
\begin{enumerate}
\item
Assume $m=(n-1)a+1$ for $a\geq 1$.
Let $z\in {\cal W}_m$ and $l=(m-1)k+1$ for some $k\geq 1$.
Define $\hat{z}\in {\cal W}_l$ by
%
%
\begin{equation}
\label{eqn:hatzsum}
\hat{z}:=\sum_{r=0}^{k-1}\sum_{i=1}^{m-1}z_m^rz_ie_{(m-1)r+i}+
z_m^ke_l
\end{equation}
where $\{e_i\}$ denotes the standard basis of ${\Bbb C}^l$.
Then $\omega_z'=\omega'_{\hat{z}}$.
\item
Let ${\cal S}_{m,n}$ denote the set of all geometric progression
states on $\con$ parametrized by $z\in {\cal W}_m$, that is,
 ${\cal S}_{m,n}:=\{\omega_z':z\in {\cal W}_m\}$.
If $l\geq m\geq 2$ satisfy that $m-1$ is a divisor of $l-1$, then 
${\cal S}_{m,n}\subset {\cal S}_{l,n}$.
\item
For any $m,l\in \{(n-1)k+1:k\geq 2\}$,
let $p:=(m-1)(l-1)+1$ and  $z\in {\cal W}_m$ and $y\in {\cal W}_l$.
Then 
$\omega_z'\sim \omega_y'$
if and only if the following equation of vectors in ${\Bbb C}^p$ holds:
%
%
\begin{equation}
\label{eqn:hatyhatz}
\sum_{r=0}^{l-2}\sum_{i=1}^{m-1}z_m^rz_ie_{(m-1)r+i}
+z_m^{l-1}e_{p}
=\sum_{r'=0}^{m-2}\sum_{i'=1}^{l-1}y_l^{r'}y_{i'}e_{(l-1)r'+i'}
+y_l^{m-1}e_{p}
\end{equation}
where $\{e_i\}$ denotes the standard basis of ${\Bbb C}^p$.
%
\item
Assume $m=(n-1)k+1$ and $k\geq 2$.
Let $z\in {\cal W}_m$
and $y=(y_1,\ldots,y_n)\in ({\Bbb C}^n)_1$.
Let $\omega_y$ be as in (\ref{eqn:first}).
Then $\omega_z'\sim \omega_y$ 
if and only if
$|y_n|<1$ and $z=\hat{y}$
where $\hat{y}\in ({\Bbb C}^m)_1$ is defined as 
%
%
\begin{equation}
\label{eqn:tildey}
\hat{y}:=\sum_{r=0}^{k-1}\sum_{j=1}^{n-1}y_n^{r}y_je_{(n-1)r+j}
+y^k_ne_m
\end{equation}
where $\{e_j\}$ denotes the standard basis of ${\Bbb C}^m$.
\end{enumerate}
\end{Thm}


%
%
\sssft{Finite correlation}
\label{subsubsection:firstthreethree}
For $2\leq n<\infty$, define
%
%
\begin{equation}
\label{eqn:idefinition}
{\cal I}_n:=\bigcup_{l\geq 0}\{1,\ldots,n\}^l,\quad
\{1,\ldots,n\}^0:=\{\emptyset\}.
\end{equation}
For $n=\infty$,
replace $\{1,\ldots,n\}$ by ${\Bbb N}:=\{1,2,\ldots\}$ in (\ref{eqn:idefinition}).
For $J=(j_1,\ldots,j_r)\in {\cal I}_n$,
we write $s_J:=s_{j_1}\cdots s_{j_r}$ 
and let $s_{\emptyset}:=I$  for convenience.
For a state $\omega$ on $\con$,
$\omega$ is said to be {\it finitely correlated} \cite{BJP}
if the dimension of $
{\cal K}(\omega):
={\rm Lin}\langle\{\pi(s_J)^*\Omega:J\in {\cal I}_n\}\rangle$
is finite
where $({\cal H},\pi,\Omega)$ denotes
the GNS representation by $\omega$.
%
%
\begin{Thm}
\label{Thm:finite}
For $n<\infty$,
any geometric progression state $\omega$ on $\con$ of order $k<\infty$ satisfies
$\dim {\cal K}(\omega)\leq k$.
Especially,  $\omega$ is finitely correlated.
\end{Thm}

\noindent
It is known that any sub-Cuntz state on $\con$ ($n<\infty$) is finitely correlated
(see Lemma \ref{lem:gpeight}(i)).
Furthermore,
it is easily shown that 
$\omega$ is a Cuntz state if and only if $\dim {\cal K}(\omega)=1$.
%
%
\ssft{Properties of state parametrization}
\label{subsection:firstfour}
In this subsection,
we show properties of the state parametrization
%
%
\begin{equation}
\label{eqn:parametrization}
{\cal W}_m\ni z\longmapsto \omega_z'\in {\cal P}(\con)
\end{equation}
where ${\cal P}(\con)$ denotes the set of all pure states on $\con$
and ${\cal W}_{\infty}:=\ltno$ when $m=\infty$ for convenience.
We consider the naturality and relevance of the parametrization
(\ref{eqn:parametrization}).
From Theorem \ref{Thm:mainb},
(\ref{eqn:parametrization})
is injective and it can be defined into
the set of all unitary equivalence classes of pure states (or irreducible representations)
of $\con$ (= the spectrum of $\con$ \cite{Dixmier})
from Theorem \ref{Thm:mainthree}.
In addition,  we show the following two properties:
(i) (\ref{eqn:parametrization}) is covariant with respect to
certain $U(n-1)$-actions.
(ii) (\ref{eqn:parametrization}) is an isomorphism of two inductive systems.
%
%
\sssft{$U(n-1)$-covariance of state parametrization}
\label{subsubsection:firstfourone}
%
We introduce two actions of unitary groups on ${\cal W}_m$
and ${\cal P}(\con)$ as follows.
For $2\leq m<\infty$,
we write the standard action of $U(m)$ on the vector space ${\Bbb C}^m$
as $gz$ for $g\in U(m)$ and $z\in {\Bbb C}^m$.
Identify $U(m-1)$ with a subgroup of $U(m)$ with respect to
the embedding 
%
%
\begin{equation}
\label{eqn:unminus}
U(m-1)\ni g\mapsto 
\left[
\begin{array}{cc}
g &0\\
0& 1\\
\end{array}
\right]\in U(m).
\end{equation}
The subgroup $U(m-1)$ of $U(m)$ also acts on ${\Bbb C}^m$.
From (\ref{eqn:unminus}) and the definition of ${\cal W}_m$,
the subset ${\cal W}_m\subset {\Bbb C}^m$
is invariant under the action of $U(m-1)$,
that is, $g{\cal W}_m\subset {\cal W}_m$ for all $g\in U(m-1)$. 

On the other hand,
let $\alpha$ denote the standard $U(n)$-action on $\con$
defined by
$\alpha_g(s_i):=\sum_{j=1}^{n}g_{ji}s_j$ for $i=1,\ldots,n$
and $g=(g_{ij})\in U(n)$.
For a state $\omega$ on $\con$,
define $\alpha_g^*(\omega):=\omega\circ \alpha_{g^{-1}}$
for $g\in U(n)$.
Identify $U(n-1)$ with a subgroup of $U(n)$ with respect to
the embedding in (\ref{eqn:unminus}) by replacing $m$ with $n$.
By this identification,
$U(n-1)$ also acts on $\con$ such that $\alpha_g(s_n)=s_n$
for any $g\in U(n-1)$.

Then we have the following result.
%
%
\begin{Thm}($U(n-1)$-covariance)
\label{Thm:cov}
\begin{enumerate}
\item
For $m=(n-1)k+1$ $(k\geq 2)$ and $z\in {\cal W}_{m}$,
the following holds:
%
%
\begin{equation}
\label{eqn:alphagstar}
\alpha_{g}^*(\omega_z') =\omega_{\tilde{g}z}'
\quad(g\in U(n-1))
\end{equation}
where $\tilde{g}=(\tilde{g}_{i,j})\in U(m-1)\subset U(m)$ is defined as
$\tilde{g}_{(n-1)a+i,(n-1)b+j}:=\delta_{ab}\,g_{ij}$ for $i,j=1,\ldots,n-1$
and $0\leq a,b\leq k-1$
and $\tilde{g}_{m,m}:=1$.
\item
For any $z\in \ltno$, the following holds:
%
%
\begin{equation}
\label{eqn:alphagstari}
\alpha_{g}^*(\omega_z') =\omega_{\tilde{g}z}'
\quad(g\in U(n-1))
\end{equation}
where $\tilde{g}=(\tilde{g}_{i,j})$
is the unitary operator on $\ltn$ defined as
$\tilde{g}_{(n-1)a+i,(n-1)b+j}:=\delta_{ab}\,g_{ij}$ for $i,j=1,\ldots,n-1$
and $a,b\geq 0$
where $\tilde{g}_{i,j}$'s  denote matrix components  of $\tilde{g}$
with respect to the standard basis of $\ltn$.
\end{enumerate}
\end{Thm}
%
%
\pr
(i)
By definition, $\tilde{g}z\in {\cal W}_m$. 
Hence $\omega'_{\tilde{g}z}$ is uniquely defined
from Theorem \ref{Thm:mainb}.
Identify $\co{m}$ with $f(\co{m})$ for $f$ in (\ref{eqn:gencc}).
Then we can verify
$\{\alpha_g^*(\omega_z')\}(t_{(n-1)r+i})=\overline{(\tilde{g}z)_{(n-1)r+i}}$
for all $r=0,1,\ldots,k-1$ and $i=1,\ldots,n-1$,
and $\{\alpha_g^*(\omega_z')\}(t_m)=\overline{(\tilde{g}z)_m}$.
Hence $\alpha_g^*(\omega_z')=\omega'_{\tilde{g}z}$
by the uniqueness of $\omega'_{\tilde{g}z}$.

\noindent
(ii) As the same token with the proof of (i),
the statement can be proved.
\qedh

\noindent
From Theorem \ref{Thm:cov}(i) ({\it resp.} (ii)),
we see that 
the parametrization (\ref{eqn:parametrization}) 
is covariant with respect to two actions of $U(n-1)$
on ${\cal W}_m$ ({\it resp.} $\ltno$) and ${\cal P}(\con)$.
When $m=n$ in Theorem \ref{Thm:cov},
$\omega_z'$ is just the Cuntz state $\omega_z$.
In this case,
it is known that
$\alpha_g^*(\omega_z)=\omega_{gz}$ for all $g\in U(n)$
and $z\in ({\Bbb C}^n)_1$ (\cite{GP08}, (1.14)).
Hence Theorem \ref{Thm:cov} is a natural generalization of 
this covariance.
%
%
\sssft{State parametrization as an isomorphism 
between two inductive systems}
\label{subsubsection:firstfourtwo}
%
For a directed set $(D,\leq)$,
a data $\{(A_{d},\varphi_{e,d}):d,e\in D\}$ is a (set-theoretical) 
{\it inductive system (or directed system)} over $(D,\leq)$
\cite{Bourbaki}
if maps
$\varphi_{e,d}:A_{d}\to A_{e}$ $(d\leq e)$
satisfy 
$\varphi_{g,e}\circ \varphi_{e,d}=\varphi_{g,d}$ $(d\leq e\leq g)$ and
$\varphi_{d,d}=id_{A_{d}}$.
For $l,k\in {\Bbb N}:=\{1,2,\ldots\}$,
if $k$ is a divisor of $l$, then we write $k\prec  l$.
We introduce two inductive systems
over the directed set $({\Bbb N},\prec)$.
%
%
\begin{Thm}
\label{Thm:inductive}
Fix $2\leq n<\infty$.
For ${\cal S}_{m,n}$ in Theorem \ref{Thm:mequivalence}(ii) and ${\cal W}_{m}$
in Theorem \ref{Thm:mainthree}(i),
define
%
%
\begin{equation}
\label{eqn:swk}
{\cal S}(k):={\cal S}_{(n-1)k+1,n},\quad 
{\cal W}(k):={\cal W}_{(n-1)k+1}\quad (k\geq 1).
\end{equation}
\begin{enumerate}
\item
$\{({\cal S}(k),\subset):k\in {\Bbb N}\}$ is an inductive system
over $({\Bbb N},\prec)$
with respect to 
inclusions $\subset$ in Theorem \ref{Thm:mequivalence}(ii).
\item
When $l\succ k$, define the map 
$\psi_{l,k}:{\cal W}(k)\to {\cal W}(l)$ by
%
%
\begin{equation}
\label{eqn:psilk}
\psi_{l,k}(z):=\hat{z}
\end{equation}
where $\hat{z}$ is as in Theorem \ref{Thm:mequivalence}(i).
Then 
$\{({\cal W}(k),\psi_{l,k}):k,l\in {\Bbb N}\}$ is
an inductive system over $({\Bbb N},\prec)$.
\item
The state parametrization 
%
%
\begin{equation}
\label{eqn:phikw}
\Phi_k:{\cal W}(k)\ni z\longmapsto \omega_z'\in {\cal S}(k)
\quad(k\geq 1)
\end{equation}
gives an isomorphism $\{\Phi_k:k\geq 1\}$
between two inductive systems in (i) and (ii).
\end{enumerate}
\end{Thm}
%
%
\pr
(i)
Assume $k\prec l$.
From Theorem \ref{Thm:mequivalence}(ii), ${\cal S}(k)\subset {\cal S}(l)$.
We see that its inclusion map
$\phi_{l,k}:{\cal S}(k)\hookrightarrow {\cal S}(l)$ 
satisfies $\phi_{m,l}\circ \phi_{l,k}=\phi_{m,k}$
when $m\succ l\succ k$.

\noindent
(ii)
By definition,
we can prove  $\psi_{m,l}\circ \psi_{l,k}=\psi_{m,k}$
when $m\succ l\succ k$.
Hence the statement holds.

\noindent
(iii)
By definition, $\Phi_k$ is a bijection.
From proofs of (i) and (ii), the statement holds.
\qedh

The paper is organized as follows.
In $\S$ \ref{section:second},
we will review known results and prepare lemmas to prove main theorems.
In $\S$ \ref{section:third},
we will prove main theorems.
In $\S$ \ref{section:fourth},
we will show examples.

%
%
\sftt{Preparations}
\label{section:second}
%
%
%
\ssft{Equivalent conditions of geometric progression state}
\label{subsection:secondone}
For convenience, we show equivalent conditions of geometric progression state here.
Let $s_1,\ldots,s_n$ denote Cuntz generators of $\con$.
From Lemma \ref{lem:homo}(ii) and Definition \ref{defi:lexico},
the following holds.
%
%
\begin{cor}
\label{cor:homoone}
Let $\omega$ be a state on $\con$ with the GNS representation
$({\cal H},\pi,\Omega)$.
\begin{enumerate}
\item
Assume $m=(n-1)k+1$ for $k\geq 1$.
For $z\in ({\Bbb C}^m)_1$,
the following are equivalent:
\begin{enumerate}
\item
$\omega$ is a geometric progression state by $z$.
\item
$\disp{\Bigl\{\sum_{r=0}^{k-1}\sum_{i=1}^{n-1}z_{(n-1)r+i}\, 
\pi(s_n^rs_i)
+z_m\pi(s_n^k)\Bigr\}\Omega=\Omega}$.
\item
$\pi(s_n^rs_i)^*\Omega=z_{(n-1)r+i}\Omega$
for 
$r=0,1,\ldots,k-1$ and 
$i=1,\ldots,n-1$, and 
$\pi(s_n^k)^*\Omega=z_{m}\Omega$.
\end{enumerate}
\item
For $z\in\ltno$, the following are equivalent:
\begin{enumerate}
\item
$\omega$ is a geometric progression state by $z$.
\item
$\disp{\sum_{r\geq 0}\sum_{i=1}^{n-1}z_{(n-1)r+i}\, 
\pi(s_n^rs_i)\Omega=\Omega}$.
\item
$\pi(s_n^rs_i)^*\Omega=z_{(n-1)r+i}\Omega$ 
for $r\geq 0$ and $i=1,\ldots,n-1$.
\end{enumerate}
\end{enumerate}
\end{cor}

\noindent
For the case of $k=1$ in Corollary \ref{cor:homoone}(i),
we obtain the equivalent conditions of Cuntz state,
that is,  the following are equivalent:
\begin{enumerate}
\item
$\omega$ is the Cuntz state by $z$.
\item
$\{z_{1}\pi(s_1)+\cdots +z_{n}\pi(s_n) \}\Omega=\Omega$.
\item
$\pi(s_i)^*\Omega=z_i\Omega$ for $i=1,\ldots,n$.
\end{enumerate}

We show relations between Cuntz generators of $\co{m}$ and $\con$.
%
%
\begin{lem}
\label{lem:words}
Let ${\cal I}_n$ be as in (\ref{eqn:idefinition}).
\begin{enumerate}
\item
Assume $m=(n-1)k+1$ for $k\geq 2$.
For $f$ in (\ref{eqn:gencc}),
we write $f(t_i)$ as $t_i$  for short.
For any $J\in {\cal I}_n$,
there exists a unique pair $(\hat{J},a)\in {\cal I}_m\times \{0,1,\ldots,k-1\}$
such that $s_J=t_{\hat{J}}s_n^a$.
\item
For $f$ in (\ref{eqn:gencci}), we write $f(t_i)$ as $t_i$ for short.
For any $J\in {\cal I}_n$,
there exists a unique pair $(\hat{J},a)\in {\cal I}_{\infty}\times {\Bbb Z}_{\geq 0}$
such that $s_J=t_{\hat{J}}s_n^a$.
\end{enumerate}
\end{lem}
%
%
\pr
See Appendix \ref{section:appone}.
\qedh

%
%
\ssft{Sub-Cuntz states}
\label{subsection:secondtwo}
In this subsection,
we review sub-Cuntz state \cite{GP08}.
For $m\geq 1$,
let ${\cal V}_{n,m}$ denote
the complex Hilbert space with the orthonormal basis
$\{e_J: J\in\{1,\ldots,n\}^m\}$,
that is, ${\cal V}_{n,m}=\ell^2(\{1,\ldots,n\}^m)\cong {\Bbb C}^{n^m}$.
Let  $({\cal V}_{n,m})_1
:=\{z\in {\cal V}_{n,m}:\|z\|=1\}$.
%
%
\begin{defi}
\label{defi:subcuntz}
For $z=\sum z_Je_J\in ({\cal V}_{n,m})_1$,
$\omega$ is a {\it sub-Cuntz state} on $\con$ by $z$ 
if $\omega$ is a state on $\con$ which satisfies the following equations:
%
%
\begin{equation}
\label{eqn:subeqn}
\omega(s_{J})=\overline{z_{J}}
\quad \mbox{for all }J\in \{1,\ldots,n\}^m
\end{equation}
where $s_J:=s_{j_1}\cdots s_{j_m}$ when
$J=(j_1,\ldots,j_m)$,
and $\overline{z_J}$ denotes the complex conjugate of $z_J$.
In this case,
$\omega$ is called a {\it sub-Cuntz state of order} $m$.
\end{defi}

\noindent
A sub-Cuntz state $\omega$ of order $1$ is just a Cuntz state.
A sub-Cuntz state $\omega$ of order $m$
is an $f$-sub-Cuntz state with respect to
the embedding $f\in {\rm Hom}(\co{n^m},\con)$
defined by
%
%
\begin{equation}
\label{eqn:ftisj}
f(t_i):=s_{J(i)}\quad(i=1,\ldots,n^m)
\end{equation}
where $J(i)=(j_1,\ldots,j_m)\in \{1,\ldots,n\}^m$ is defined 
as $i=\sum_{r=1}^{m}(j_r-1) n^{m-r}+1$.
For example,
if $(n,m)=(2,2)$,
then we obtain $(J(1),J(2),J(3),J(4))=((1,1),(1,2),(2,1),(2,2))$.

We identify 
${\cal V}_{n,m}$ with $({\cal V}_{n,1})^{\otimes m}$
by the correspondence between bases
$e_{J}\mapsto e_{j_1}\otimes \cdots \otimes e_{j_m}$
for $J=(j_1,\ldots,j_m)\in\{1,\ldots,n\}^m$.
From this identification,
we obtain
${\cal V}_{n,m}\otimes {\cal V}_{n,l}={\cal V}_{n,m+l}$
 for any $m,l\geq 1$.
Then the following hold.
%
%
\begin{Thm}
\label{Thm:maintwo}
\begin{enumerate}
\item
(\cite{GP08}, Fact 1.3)
For any $z\in ({\cal V}_{n,m})_1$,
a sub-Cuntz state on $\con$ by $z$ exists.
\item
(\cite{GP08}, Theorem 1.4)
For a sub-Cuntz state $\omega$ on $\con$ by $z\in ({\cal V}_{n,m})_1$,
$\omega$ is unique if and only if  $z$ is {\it nonperiodic},
that is, $z=x^{\otimes p}$ for some $x$ implies $p=1$.
In this case, $\omega$ is pure and we write it as $\tilde{\omega}_z$.
\item
(\cite{GP08}, Theorem 1.5)
Let $p\geq 2$ and $z=x^{\otimes p}$
for a nonperiodic element $x\in ({\cal V}_{n,m'})_1$.
If $\omega$ is a sub-Cuntz state on $\con$ by $z$,
then $\omega$ is a convex hull of 
sub-Cuntz states by $e^{2\pi j\sqrt{-1}/p}x$ for $j=1,\ldots,p$.
\item
(\cite{GP08}, Theorem 1.7)
For $z,y\in \bigcup_{m\geq 1}({\cal V}_{n,m})_1$,
assume that both $z$ and $y$ are nonperiodic.
Then the following are equivalent:
\begin{enumerate}
\item
GNS representations by $\tilde{\omega}_z$ and $\tilde{\omega}_y$
are unitarily equivalent. 
\item
$z$ and $y$ are {\it conjugate},
that is, 
$z=y$, or 
$z=x_1\otimes x_2$ and $y=x_2\otimes x_1$
for some $x_1, x_2\in \bigcup_{m\geq 1}({\cal V}_{n,m})_1$.
\end{enumerate}
\end{enumerate}
\end{Thm}

\noindent
About concrete examples, see Example \ref{ex:first}
(see also $\S$ 4 in \cite{GP08}).

%
%
\begin{lem}
\label{lem:gpeight}
\begin{enumerate}
\item
(\cite{GP08}, Lemma 2.4(i))
When $n<\infty$,
any sub-Cuntz state on $\con$ is finitely correlated.
\item
(\cite{GP08}, Lemma 2.4(ii))
If $\omega$ is a sub-Cuntz state with the GNS representation
$({\cal H},\pi,\Omega)$,
then ${\rm Lin}\langle\{\pi(s_J)\Omega:J\in {\cal I}_n\}\rangle$
is dense in ${\cal H}$.
\item
(\cite{GP08}, Theorem 2.3)
Fix $m\geq 1$.
Let $\omega$ be a state on $\con$ with the GNS representation
$ ({\cal H}, \pi,\Omega)$.
For $z=\sum z_Je_J\in ({\cal V}_{n,m})_1$,
the following are equivalent:
\begin{enumerate}
\item
$\omega$ is a sub-Cuntz state by $z$.
\item
$\Omega=\pi(s(z))\Omega$
where $s(z):=\sum z_Js_J$.
\item
$\pi(s_{J})^*\Omega=z_{J}\Omega$
for all $J\in \{1,\ldots,n\}^m$.
\end{enumerate}
\end{enumerate}
\end{lem}


%
%
\ssft{GNS representations by geometric progression states}
\label{subsection:secondthree}
In this subsection,
we show properties of GNS representations by geometric progression states.
Let ${\cal I}_n$ be as in (\ref{eqn:idefinition}).
For $J=(j_1,\ldots,j_r)\in {\cal I}_n$ and $z=(z_1,\ldots,z_n)\in {\Bbb C}^n$,
we write $z_J:=z_{j_1}\cdots z_{j_r}$ and 
$z_{\emptyset}:=1$ for convenience.
%
%
\begin{lem}
\label{lem:reductionthree}
Assume $m=(n-1)k+1$ for $k\geq 2$.
For $z\in ({\Bbb C}^m)_1$,
let $\omega$ be a geometric progression state on $\con$ by $z$
with the GNS representation $({\cal H},\pi,\Omega)$.
\begin{enumerate}
\item
For any $J\in {\cal I}_n$,
$\pi(s_J)^*\Omega\in 
{\rm Lin}\langle\{\pi(s_n^a)^*\Omega:a=0,\ldots,k-1\}\rangle $.
\item
For any $J\in {\cal I}_n$,
 $\pi(s_J)^*\Omega\in 
{\rm Lin}\langle\{\pi(s_L)\Omega:L\in {\cal I}_n\}\rangle $.
\item
${\rm Lin}\langle\{\pi(s_L)\Omega:L\in {\cal I}_n\}\rangle$ is dense in ${\cal H}$. 
\end{enumerate}
\end{lem}
%
%
\pr
(i)
From Lemma \ref{lem:words}(i) and 
Lemma \ref{lem:homo}(ii)((a)$\Rightarrow$(c)),
we obtain
%
%
\begin{equation}
\label{eqn:pijo}
\pi(s_J)^*\Omega
=\pi(t_{\hat{J}}s_n^a)^*\Omega
=
\pi(s_n^a)^*\pi(t_{\hat{J}})^*\Omega
=\overline{z_{\hat{J}}} \pi(s_n^a)^*\Omega
\end{equation}
for some $(\hat{J},a)$ with $0\leq a\leq k-1$.
Hence the statement holds.

\noindent
(ii)
From Corollary \ref{cor:homoone}(i)((a)$\Rightarrow$(b)), we can prove
%
%
\begin{equation}
\label{eqn:pisna}
\pi(s_n^a)^*\Omega
=\Bigl\{\sum_{r=a}^{k-1}\sum_{i=1}^{n-1}z_{(n-1)r+i}\pi(s_n^{r-a}s_i)
+z_m \pi(s_n^{k-a})\Bigr\}\Omega.
\end{equation}
From this and (i), the statement holds.

\noindent
(iii)
Since ${\rm Lin}\langle\{\pi(s_Js_K^*)\Omega: J,K\in {\cal I}_n\}\rangle $
is dense in ${\cal H}$,
the statement holds from (ii).
\qedh

%
%
\begin{lem}
\label{lem:reductionthreei}
For $z\in\ltno$,
let $\omega$ be a geometric progression state on $\con$ by $z$
with the GNS representation $({\cal H},\pi,\Omega)$.
\begin{enumerate}
\item
For any $J\in {\cal I}_n$,
$\pi(s_J)^*\Omega\in 
{\rm Lin}\langle\{\pi(s_n^a)^*\Omega:a\geq 0\}\rangle $.
\item
For any $J\in {\cal I}_n$,
 $\pi(s_J)^*\Omega\in
\overline{{\rm Lin}\langle\{\pi(s_L)\Omega:L\in {\cal I}_n\}\rangle }$.
\item
${\rm Lin}\langle\{\pi(s_L)\Omega:L\in {\cal I}_n\}\rangle$
is dense in ${\cal H}$. 
\end{enumerate}
\end{lem}
%
%
\pr
By replacing ``Lemma \ref{lem:words}(i)"
in the proof of Lemma \ref{lem:reductionthree}
with ``Lemma \ref{lem:words}(ii)",
all statements can be verified.
\qedh

%
%
\ssft{General properties of $f$-sub-Cuntz states}
\label{subsection:secondfour}
%
%
\begin{lem}
\label{lem:equivfive}
Assume that $A$ is a unital C$^*$-algebra
and $B$ is a unital C$^*$-subalgebra of $A$.
For two states $\omega$ and $\omega'$ on $A$,
if $\omega$ is pure and 
the restriction $\pi_{\omega}|_B$ is irreducible,
then $\omega\sim\omega'$ implies $\omega|_B\sim\omega'|_B$.
\end{lem}
%
%
\pr
Assume $\omega\sim\omega'$.
Since $\omega$ is pure, $\omega'$ is also pure.
Hence there exists an irreducible representation $\pi$
of $A$ on a Hilbert space ${\cal H}$ with
two cyclic unit vectors $\Omega$ and $\Omega'$
such that 
$\omega=\langle\Omega|\pi(\cdot)\Omega\rangle$
and 
$\omega'=\langle\Omega'|\pi(\cdot)\Omega'\rangle$.
By assumption,
$\overline{\pi(B)v}={\cal H}$ for any nonzero vector $v\in {\cal H}$.
Hence 
$\overline{\pi(B)\Omega}
={\cal H}
=\overline{\pi(B)\Omega'}$.
From this,
$(\overline{\pi(B)\Omega},\pi|_B)$
and 
$(\overline{\pi(B)\Omega'},\pi|_B)$ are unitarily equivalent
as representations of $B$.
Since
$\omega|_B=\langle\Omega|\pi|_B(\cdot)\Omega \rangle$
and 
$\omega'|_B=\langle\Omega'|\pi|_B(\cdot)\Omega' \rangle$,
we obtain $\omega|_B\sim\omega'|_B$.
\qedh

\noindent
Lemma \ref{lem:equivfive} does not hold
when $\pi_{\omega}|_B$ is not irreducible.
For example, let $A:=M_{3}({\Bbb C})\curvearrowright{\Bbb C}^3$ and 
$B:={\Bbb C}\oplus M_2({\Bbb C})\subset A$,
and let $e_1,e_2,e_3$ denote
the standard basis of ${\Bbb C}^3$.
Define two states $\omega:=\langle e_1|(\cdot)e_1\rangle $
and 
$\omega':=\langle e_2|(\cdot)e_2\rangle $ on $A$.
Then $\omega\sim \omega'$, but
$\omega|_{B}\not \sim \omega'|_{B}$.

Let $V^{(m)}:={\Bbb C}^m$ when $m<\infty$
and $V^{(\infty)}=\ltn$, and define
$V^{(m)}_1:=\{v\in V^{(m)}:\|v\|=1\}$
for $2\leq m\leq \infty$.
Let $({\cal H}_{\omega},\pi_{\omega},\Omega_{\omega})$
denote the GNS representation by a state $\omega$.
%
%
\begin{cor}
\label{cor:rest}
For $2\leq m\leq \infty$, 
fix $f\in {\rm Hom}(\co{m},\con)$ and identify $\co{m}$ with $f(\co{m})$.
For $z,y\in V_1^{(m)}$,
let $\omega$ and $\omega'$ be $f$-sub-Cuntz states
by $z$ and $y$, respectively such that 
$\pi_{\omega}|_{\co{m}}$ is irreducible.
If $\omega\sim \omega'$, then $z= y$.
\end{cor}
%
%
\pr
Assume $\omega\sim \omega'$.
Remark that restrictions $\omega|_{\co{m}}$ 
and $\omega'|_{\co{m}}$ are Cuntz states $\omega_z$ 
and $\omega_y$ on $\co{m}$, respectively by definition.
From this and Lemma \ref{lem:equivfive},
$\omega_z=\omega|_{\co{m}}\sim \omega'|_{\co{m}}=\omega_y$.
This is equivalent to $z=y$.
\qedh

For given two embeddings $f$ and $g$,
we show a sufficient condition of the equivalence
between $f$-sub-Cuntz states and $g$-sub-Cuntz states.
%
%
\begin{Thm}
\label{Thm:gequiv}
For $2\leq m,l\leq \infty$,
let $f\in {\rm Hom}(\co{m},\con)$
and $g\in {\rm Hom}(\co{l},\con)$ and 
let $t_1,\ldots,t_m$ and $u_1,\ldots,u_l$
denote Cuntz generators of $\co{m}$ and $\co{l}$,
respectively.
For $z\in V^{(m)}_1$ and $y\in V^{(l)}_1$,
assume 
%
%
\begin{equation}
\label{eqn:ftzguy}
f(t(z))=g(u(y))
\end{equation}
for
$t(z):=\sum_{j=1}^{m}z_jt_j\in \co{m}$
and 
$u(y):=\sum_{i=1}^{l}y_iu_i\in \co{l}$.
\begin{enumerate}
\item
For a state $\omega$ on $\con$,
$\omega$ is an $f$-sub-Cuntz state by $z$
if and only if 
$\omega$ is a $g$-sub-Cuntz state by $y$.
\item
If an $f$-sub-Cuntz state by $z$
or a $g$-sub-Cuntz state by $y$
is unique,
then they coincide as a state on $\con$.
Especially, they are equivalent.
\end{enumerate}
\end{Thm}
%
%
\pr
(i)
From Lemma \ref{lem:homo}(ii)((a)$\Leftrightarrow$(b))
and (\ref{eqn:ftzguy}),
the statement holds.

\noindent
(ii)
From (i), 
if one is unique, then so is other.
Hence the statement holds from (i).
\qedh

%
%
\sftt{Proofs of main theorems}
\label{section:third}
%
%
%
\ssft{Proofs of Theorem \ref{Thm:mainb}\, and Theorem \ref{Thm:mixture}}
\label{subsection:thirdone}
In order to prove Theorem \ref{Thm:mainb},
we show formulas of explicit values of \gprs s.
Let ${\cal I}_n$ be as in (\ref{eqn:idefinition}).
%
%
\begin{Thm}
\label{Thm:explicit}
Fix $k\geq 2$ and let $m=(n-1)k+1$.
For $z=(z_1,\ldots,z_m)\in ({\Bbb C}^m)_1$,
let $\omega$ be a \gprs\ on $\con$ by $z$.
If $|z_m|<1$,
then the following holds:
\begin{enumerate}
\item
For $f$ in (\ref{eqn:gencc}), we write $f(t_i)$ as $t_i$ for short.
For $J,K\in {\cal I}_n$,
there exist unique 
$(\hat{J},a)$ and $(\hat{K},b)$ 
in ${\cal I}_m\times \{0,1,\ldots,k-1\}$
such that 
$s_J=t_{\hat{J}}s_n^a$, $s_K=t_{\hat{K}}s_n^b$
and 
$\omega(s_Js_K^*)=
\overline{z_{\hat{J}}}z_{\hat{K}}\,\omega(s_n^a(s_n^b)^*)$.
\item
When $0\leq a\leq b\leq k-1$,
%
%
\begin{equation}
\label{eqn:omsum}
\omega(s_n^a(s_n^b)^*)=
\sum_{j=1}^{(n-1)(k-b)}\overline{z_{(n-1)a+j}}\,z_{(n-1)b+j}
+\frac{|z_m|^2Z_{b-a}+\overline{z_m Z_{k-b+a}}}{1-|z_m|^2}
\end{equation}
where $Z_c\in {\Bbb C}$ $(0\leq c\leq k)$ is defined as
%
%
\begin{equation}
\label{eqn:zc}
Z_c:=\sum_{r=1}^{(n-1)(k-c)}\overline{z_{(n-1)c+r}}z_r
\quad (0\leq c\leq k-1),\quad
Z_k:=0.
\end{equation}
\end{enumerate}
\end{Thm}
%
%
\pr
(i)
From Lemma \ref{lem:words}(i) 
and Lemma \ref{lem:homo}(ii)((a)$\Rightarrow$(c)), the statement holds.

\noindent
(ii)
Let $\Theta_{a,b}:=\omega(s_n^a(s_n^b)^*)$.
Assume $k-2\geq b\geq a\geq 0$.
By Corollary \ref{cor:homoone}(i)((a)$\Rightarrow$(c)) 
and $\sum_{i=1}^{n}s_is_i^*=I$,
we obtain
%
%
\begin{equation}
\label{eqn:thetaab}
\Theta_{a,b}
=\sum_{i=1}^{n}
\omega(s_n^as_is_i^*(s_n^b)^*)
=\sum_{i=1}^{n-1}
\overline{z_{(n-1)a+i}}z_{(n-1)b+i}
+\Theta_{a+1,b+1}.
\end{equation}
By repetition to $\Theta_{a+1,b+1}$, we obtain
%
%
\begin{equation}
\label{eqn:thetabb}
\begin{array}{rl}
\Theta_{a,b}
=&
\disp{\sum_{r=a}^{k-1-b+a}\sum_{i=1}^{n-1}
\overline{z_{(n-1)r+i}}z_{(n-1)(r+b-a)+i}
+\Theta_{k-b+a,k}}\\
\\
=&
\disp{\sum_{j=(n-1)a+1}^{(n-1)(k-b+a)}
\overline{z_j}z_{(n-1)(b-a)+j}
+z_m\omega(s_n^{k-b+a}).}\\
\end{array}
\end{equation}
On the other hand,
from Corollary \ref{cor:homoone}(i)((a)$\Rightarrow$(b)),
%
%
\begin{equation}
\label{eqn:omegasna}
\begin{array}{rl}
\omega(s_n^a)
=&\disp{
\sum_{i=1}^{n-1}
\sum_{j=a}^{k-1}
\overline{z_{(n-1)j+i}}
\omega(t_{(n-1)(j-a)+i}^*)
+
\overline{z_m}\omega((s_n^{k-a})^*)
}\\
=&\disp{
Z_a+\overline{z_m\omega(s_n^{k-a})}}\\
\end{array}
\end{equation}
where  we use $t_{(n-1)j+i}^*s_n^a=t_{(n-1)(j-a)+i}^*$ when $j\geq a$
and $t_{(n-1)j+i}^*s_n^a=0$ when $j< a$.
By replacing $a$ with $k-a$ in (\ref{eqn:omegasna}),
we obtain
$\omega(s_n^{k-a})=Z_{k-a}+\overline{z_m\omega(s_n^{a})}$.
By substituting this into (\ref{eqn:omegasna}), we obtain
%
%
\begin{equation}
\label{eqn:omesnaz}
\omega(s_n^a)=
\frac{\overline{z_m Z_{k-a}}+Z_a}{1-|z_m|^2}\quad
(0\leq a\leq k).
\end{equation}
From this and (\ref{eqn:thetabb}),
the statement is verified.
\qedh

%
%
\begin{Thm}
\label{Thm:oniunique}
For $z\in\ltno$,
let $\omega$ be a \gprs \ on $\con$ by $z$.
For $f$ in (\ref{eqn:gencci}),
we write $f(t_i)$ as $t_i$ for short.
For $J,K\in {\cal I}_n$,
there exist unique
$(\hat{J},a)$ and $(\hat{K},b)$ in ${\cal I}_{\infty}\times
{\Bbb Z}_{\geq 0}$  such that 
$s_J=t_{\hat{J}}s_n^a$, $s_K=t_{\hat{K}}s_n^b$
and 
%
%
\begin{equation}
\label{eqn:oniuniqeu}
\omega(s_Js_K^*)=
\overline{z_{\hat{J}}}z_{\hat{K}}
\sum_{j\geq 1}\overline{z_{(n-1)a+j}}\,z_{(n-1)b+j}.
\end{equation}
\end{Thm}
%
%
\pr
Let $J,K\in {\cal I}_n$.
From Lemma \ref{lem:words}(ii) and 
Lemma \ref{lem:homo}(ii)((a)$\Rightarrow$(c)), 
$\omega(s_Js_K^*)
=\overline{z_{\hat{J}}}z_{\hat{K}}
\omega(s_n^a(s_n^b)^*)$ for some $(\hat{J},a)$
and $(\hat{K},b)$.
From Corollary \ref{cor:homoone}(ii), 
we can prove
$\omega(s_n^a(s_n^b)^*)
=\sum_{j\geq 1}\overline{z_{(n-1)a+j}}z_{(n-1)b+j}$.
Hence (\ref{eqn:oniuniqeu}) is proved.
\qedh

\noindent
{\it Proof of Theorem \ref{Thm:mixture}.}
If $|z_m|=1$, then the equation 
$z_m\pi(s_n^k)\Omega=\Omega$ holds 
from Corollary \ref{cor:homoone}(i)((a)$\Rightarrow$(b)).
Hence the state is a sub-Cuntz state 
by $z_m e_n^{\otimes k}\in ({\cal V}_{n,k})_1$ from
Lemma \ref{lem:gpeight}(iii)((b)$\Rightarrow$(a)).
From Theorem \ref{Thm:maintwo}(iii),
the statement holds. 
\qedh

\noindent
{\it Proof of Theorem \ref{Thm:mainb}.}
(i)
($\Leftarrow$)
From Theorem  \ref{Thm:explicit}, 
the statement holds.

\noindent
($\Rightarrow$)
From Theorem \ref{Thm:mixture}, the statement holds.

\noindent
(ii)
From Theorem \ref{Thm:oniunique}, the statement holds.
\qedh

%
%
\ssft{Proof of Theorem \ref{Thm:mainthree}}
\label{subsection:thirdtwo}
In order to prove Theorem \ref{Thm:mainthree},
we show the following theorem.
%
%
\begin{Thm}
\label{Thm:restriction}
For a state $\omega$, let $\pi_{\omega}$ denote 
the GNS representation by $\omega$.
\begin{enumerate}
\item
Identify $\coni$ with $f(\coni)$ for $f$ in (\ref{eqn:gencci}).
For $z\in \ltno$,
let $\omega_z'$ be as in Theorem \ref{Thm:mainb}(ii)
and let $\omega_z$ be the Cuntz state on $\coni$ by $z$.
Then the restriction $\pi_{\omega_z'}|_{\coni}$ of  $\pi_{\omega_z'}$
to $\coni$ is unitarily equivalent to $\pi_{\omega_z}$,  that is,
$\pi_{\omega_z'}|_{\coni}\sim \pi_{\omega_z}$.
Especially, $\pi_{\omega_z'}|_{\coni}$ is irreducible.
\item
Assume $m=(n-1)k+1$ for $k\geq 2$.
Identify $\co{m}$ with $f(\co{m})$ for $f$ in (\ref{eqn:gencc}).
For $z\in {\cal W}_m$,
let $\omega_z'$ be as in Theorem \ref{Thm:mainb}(i) and
let $\omega_z$ be the Cuntz state  on $\co{m}$ by $z$. 
Then the restriction $\pi_{\omega_z'}|_{\co{m}}$ of $\pi_{\omega_z'}$
to $\co{m}$
is unitarily equivalent to $\pi_{\omega_z}$,
that is,
$\pi_{\omega_z'}|_{\co{m}}\sim \pi_{\omega_z}$.
Especially, 
$\pi_{\omega_z'}|_{\co{m}}$ is irreducible.
\end{enumerate}
\end{Thm}
%
%
\pr
Let $({\cal H},\pi,\Omega)$ denote the GNS representation by $\omega_z'$
and we write $\pi(s_i)$ as $s_i$ for short.

\noindent
(i)
It is sufficient to show that $\coni\Omega$ is dense in ${\cal H}$.
From Lemma \ref{lem:reductionthreei}(iii),
it suffices to show 
$s_J\Omega\in \overline{\coni\Omega}$ for all $J$.
From Lemma \ref{lem:words}, for any $J\in {\cal I}_n$,
$s_J\Omega=t_{\hat{J}}s_n^a\Omega\in \coni s_n^a\Omega$ for some $(\hat{J},a)$.
Hence
it is enough to show $s_n^a\Omega
\in \overline{\coni\Omega}$ for any $a$.
Since $s_n^a=s_n^a\sum_{i=1}^ns_is_i^*
=\sum_{i=1}^{n-1}t_{(n-1)a+i}t_i^*+s_n^{a+1}s_n^*$,
we obtain
%
%
\begin{equation}
\label{eqn:plusone}
s_n^a\Omega
=\sum_{i=1}^{n-1}\overline{z_i}t_{(n-1)a+i}\Omega
+
s_n^{a+1}s_n^*\Omega
\end{equation}
where we use Lemma \ref{lem:homo}(ii)((a)$\Rightarrow$(c)).
Since 
$s_n^{a+1}s_n^*\Omega=
s_n^{a+1}\sum_{i=1}^ns_is_i^*s_n^*\Omega$ in
(\ref{eqn:plusone}), we obtain
%
%
\begin{equation}
\label{eqn:snaomega}
\begin{array}{rl}
s_n^a\Omega=&
\disp{\sum_{i=1}^{n-1}\overline{z(0,i)}t_{(n-1)a+i}\Omega
+
\sum_{i=1}^{n-1}\overline{z(1,i)}t_{(n-1)(a+1)+i}\Omega
+
s_n^{a+2}(s_n^2)^*\Omega
}\\
=&
\disp{\sum_{r=0}^{1}
\sum_{i=1}^{n-1}\overline{z(r,i)}t_{(n-1)(a+r)+i}\Omega
+
s_n^{a+2}(s_n^2)^*\Omega
}\\
=&\cdots \\
=&
\disp{\sum_{r=0}^{R-1}
\sum_{i=1}^{n-1}\overline{z(r,i)}t_{(n-1)(a+r)+i}\Omega
+
s_n^{a+R}(s_n^R)^*\Omega
}\\
\end{array}
\end{equation}
for all integer $R\geq 1$
where $z(r,i):=z_{(n-1)r+i}$.
Since 
$\|s_n^{a+R}(s_n^R)^*\Omega\|^2
=\|(s_n^R)^*\Omega\|^2=\sum_{r\geq 1}|z_{(n-1)R+r}|^2$
from Theorem \ref{Thm:oniunique},
$s_n^{a+R}(s_n^R)^*\Omega\to 0$ when $R\to\infty$.
From this and (\ref{eqn:snaomega}),
we obtain
$s_n^a\Omega=
\sum_{r\geq 0}
\sum_{i=1}^{n-1}\overline{z(r,i)}t_{(n-1)(a+r)+i}\Omega
\in \overline{\coni\Omega}$ for any $a\geq 1$.

\noindent
(ii)
From Lemma \ref{lem:reductionthree}(iii),
it is sufficient to show $s_J\Omega\in \co{m}\Omega$ for all $J$.
From Lemma \ref{lem:words},
for any $J\in {\cal I}_n$,
$s_J\Omega=t_{\hat{J}}s_n^a\Omega\in \co{m}s_n^a\Omega$
for some $\hat{J}$ and $0\leq a\leq  k-1$.
Hence
it is sufficient to show $s_n^a\Omega
\in \co{m}\Omega$ for any $0\leq a\leq k-1$.
By using
$s_n^a=\sum_{i=1}^{n-1}t_{(n-1)a+i}t_i^*+s_n^{a+1}s_n^*$
and the analogy of (\ref{eqn:snaomega}),
we can prove
%
%
\begin{equation}
\label{eqn:snam}
s_n^a\Omega
=
\sum_{r=0}^{k-a-1}
\sum_{i=1}^{n-1}\overline{z(r,i)}t_{(n-1)(a+r)+i}\Omega
+t_m(s_n^{k-a})^*\Omega.
\end{equation}
On the other hand,
$(s_n^{k-a})^*\Omega=
\sum_{r=k-a}^{k-1}\sum_{i=1}^{n-1}\overline{z(r,i)}
s_n^{r-(k-a)}s_i\Omega
+\overline{z_m}s_n^a\Omega$.
By substituting this into (\ref{eqn:snam}), we obtain
%
%
\begin{equation}
\label{eqn:snagg}
\begin{array}{rl}
s_n^a\Omega
=&
\disp{g
\sum_{i=1}^{n-1}
\Biggl[
\sum_{r=0}^{k-a-1}
\overline{z(r,i)}t_{(n-1)(a+r)+i}
+t_m
\sum_{r=k-a}^{k-1}\overline{z(r,i)}t_{(n-1)(r-(k-a))+i}
\Biggr]\Omega}
\end{array}
\end{equation}
where 
$g:=I+\sum_{l\geq 1}(\overline{z_m}t_m)^l\in \co{m}$.
Hence
$s_n^a\Omega\in \co{m}\Omega$ for all $0\leq a\leq k-1$.
\qedh

\noindent
{\it Proof of Theorem \ref{Thm:mainthree}.}
(i)
It is sufficient to show that 
$\omega_z'\sim \omega_y'$ implies $z= y$.
From Theorem \ref{Thm:restriction}(ii),
this holds from 
Lemma \ref{lem:equivfive} and Corollary \ref{cor:rest}.

\noindent
(ii)
As the same token with the proof of (i), the statement holds
from Theorem \ref{Thm:restriction}(i).
\qedh

%
%
\ssft{Proof of Theorem \ref{Thm:reductionequiv}}
\label{subsection:thirdthree}
Assume $m=(n-1)k+1$ for $k\geq 1$.
Let $s_1,\ldots,s_n$, $u_1,\ldots,u_m$ and $t_1,t_2\ldots$
denote Cuntz generators of $\con$,
$\co{m}$ and $\coni$, respectively.
Identify $\co{m}$ and $\coni$ with
images of geometric progression embeddings in $\con$.
\\

\noindent
{\it Proof of Theorem \ref{Thm:reductionequiv}(i).}
Let $t(z):=\sum_{i\geq 1}z_it_i$ and $s(y):=\sum_{i= 1}^ny_is_i$.

\noindent
($\Rightarrow$)
Assume  $\omega_z'\sim \omega_y$. 
Then we can assume that $\con$ acts on a Hilbert space ${\cal H}$
with two cyclic unit vectors $\Omega$ and $\Omega'$
such that $t(z)\Omega=\Omega$ and $s(y)\Omega'=\Omega'$.
Remark that $\omega_z'=\langle\Omega|(\cdot)\Omega\rangle$
and $\omega_y=\langle\Omega'|(\cdot)\Omega'\rangle$ from
Corollary \ref{cor:homoone}(ii)((b)$\Rightarrow$(a))
and Lemma \ref{lem:gpeight}(iii)((b)$\Rightarrow$(a)).
Let $X:=\langle \Omega |\Omega'\rangle$.
Then
%
%
\begin{equation}
\label{eqn:innerxy}
X=\langle t(z)\Omega|\Omega'\rangle 
=\sum_{i\geq 1}\overline{z_i}
\langle t_i\Omega|\Omega'\rangle =\alpha X
\end{equation}
where $\alpha:=
\sum_{r\geq 0}\sum_{j=1}^{n-1}
\overline{z_{(n-1)r+j}}y_{n}^ry_j\in {\Bbb C}$.
Since $\langle s_J\Omega|\Omega'\rangle =y_JX$
for all $J$,
if $X=0$, then $\omega'_z\not\sim\omega_y$.
Hence $X\ne 0$.
From this and (\ref{eqn:innerxy}), $\alpha=1$.
This implies $|y_n|<1$.
Hence $\|\tilde{y}\|=1$ and 
we can write $\alpha=\langle z|\tilde{y}\rangle$.
From $\langle z|\tilde{y}\rangle=1$ and $\|z\|=\|\tilde{y}\|=1$, 
we obtain $z=\tilde{y}$.

\noindent
($\Leftarrow$)
Assume that $|y_n|<1$ and $z=\tilde{y}$.
Let $({\cal H},\pi,\Omega')$ denote
the GNS representation by $\omega_y$.
Then $\pi(s(y))\Omega'=\Omega'$.
By assumption,
we can verify $\pi(t(z))\Omega'=\Omega'$.
From this and Corollary \ref{cor:homoone}(ii)((b)$\Rightarrow$(a)),
$\omega_z'= \omega_y$.
Especially, $\omega_z'\sim \omega_y$.
\qedh


\noindent
{\it Proof of Theorem \ref{Thm:reductionequiv}(ii).}
($\Rightarrow$)
In $\con$, we obtain
$t_{(m-1)r+i}= u_m^{r}u_i$ for $r\geq 0$ and $i=1,\ldots,m-1$.
From this,
we can verify that 
$\omega_z'|_{\co{m}}$ is the \gprs\ on $\co{m}$ by
$z$.
From Theorem \ref{Thm:restriction}(ii) and Lemma \ref{lem:equivfive},
$\omega_z'\sim \omega_y'$
implies 
%
%
\begin{equation}
\label{eqn:mequiv}
\omega_z'|_{\co{m}}\sim \omega_y'|_{\co{m}}.
\end{equation}
By definition,
 $\omega_y'|_{\co{m}}$ is the Cuntz state on $\co{m}$
by $y$.
Applying Theorem \ref{Thm:reductionequiv}(i)
to (\ref{eqn:mequiv}) by replacing $(\con,\coni)$
with $(\co{m},\coni)$,
we obtain
$z=\tilde{y}$.

\noindent
($\Leftarrow$)
Assume that $\con$ acts on a Hilbert space ${\cal H}$
with a cyclic unit vector $\Omega'$
such that $u(y)\Omega'=\Omega'$
where $u(y):=\sum_{j=1}^m y_ju_j$.
Then $\omega_y'=\langle\Omega'|(\cdot)\Omega'\rangle$.
Define $Y:=y_m$ and $y':=y-Ye_m$.
If $z=\tilde{y}$, then we obtain
%
%
\begin{equation}
\label{eqn:tzy}
t(z)=
t(\tilde{y})=\disp{
\sum_{r\geq 0}
Y^{r}u_m^ru(y')}
\end{equation}
where we use $t_{(m-1)r+j}=u_m^ru_{j}$.
Let $U:= Yu_m\in \co{m}\subset \con$.
Then we obtain
$t(z)=(I-U)^{-1}u(y')$.
From (\ref{eqn:tzy})
 and $u(y')\Omega'=(I-Yu_m)\Omega'=(I-U)\Omega'$, 
we obtain 
$t(z)\Omega'=(I-U)^{-1}u(y')\Omega'=\Omega'$.
From this,
Lemma \ref{lem:homo}(ii)((b)$\Rightarrow$(a))
and the uniqueness of $\omega_z'$,
we obtain $\omega_z'=\omega'_y$. 
Especially, 
$\omega_z'\sim \omega_y'$.
\qedh

%
%
\ssft{Proofs of Theorem \ref{Thm:mequivalence}\, and
Theorem \ref{Thm:finite}}
\label{subsection:thirdfour}

\noindent
{\it Proof of Theorem \ref{Thm:mequivalence}.}
(i)
Since
$l=(m-1)k+1=(n-1)ak+1$,
there exists the \gpre\ of $\co{l}$ into $\con$.
By assumption and the choice of $z$,
$|\hat{z}_l|=|z_m^k|<1$.
Hence $\hat{z}\in {\cal W}_l$.
Let $t_1,\ldots,t_m$ and $u_1,\ldots,u_{l}$ 
denote Cuntz generators of $\co{m}$ and $\co{l}$, respectively.
Identify $\co{m}$ and $\co{l}$ with 
images of \gpre s in $\con$, respectively.
Let $T:=\sum_{r=0}^{k-1}\sum_{i=1}^{m-1}
z_m^r z_it_m^r t_i+z_m^k t_m^k\in \co{m}\subset \con$.
Then we see $T=u(\hat{z}):=\sum_{j=1}^l\hat{z}_ju_j\in \co{l}\subset \con$.
Let $z':=z-z_me_m\in {\Bbb C}^m$ and $Z:=z_mt_m\in \co{m}\subset \con$.
Then we can rewrite
$T=\sum_{r=0}^{k-1}(z_m t_m)^r\, t(z')+(z_m t_m)^k=
(I-Z^k)(I-Z)^{-1}t(z')+Z^k$.
Let $({\cal H},\pi,\Omega)$ denote the GNS representation by $\omega_z'$.
We write $\pi(s_i)$ as $s_i$ for short.
Since $\Omega=t(z)\Omega=\{t(z')+z_mt_m\}\Omega$,
$t(z')\Omega=(I-z_mt_m)\Omega=(I-Z)\Omega$.
By using these, we obtain
$u(\hat{z})\Omega=T\Omega=\{(I-Z^k)(I-Z)^{-1}t(z')+Z^k\}\Omega
=\Omega$.
From this,
Lemma \ref{lem:homo}(ii)((b)$\Rightarrow$(a))
and the uniqueness of $\omega_{\hat{z}}'$,
the statement holds.

\noindent
(ii)
From (i), the statement holds.

\noindent
(iii)
By assumption,
$m=(n-1)a+1$
and $l=(n-1)b+1$ for some $a,b\geq 1$.
Since $p=(n-1)^2ab+1$,
we can define the \gpre\  of $\co{p}$
into $\con$.
Therefore
${\cal S}_{m,n}\cup {\cal S}_{l,n}\subset {\cal S}_{p,n}$ from (ii).
From (\ref{eqn:hatzsum}),
we obtain $\hat{z},\hat{y}\in {\cal W}_p$
such that 
$\omega_z'=\omega_{\hat{z}}'$
and 
$\omega_y'=\omega_{\hat{y}}'$.
Hence 
$\omega_z'\sim \omega_y'$
if and only if 
$\omega_{\hat{z}}'\sim \omega_{\hat{y}}'$.
This is equivalent to $\hat{z}=\hat{y}$ from 
Theorem \ref{Thm:mainthree}(i).
Since $\hat{z}$ and $\hat{y}$ coincide with the l.h.s. and the r.h.s. of   
(\ref{eqn:hatyhatz}), respectively,  the statement holds.

\noindent
(iv)
From Theorem \ref{Thm:reductionequiv}(ii),
$\omega_z'\sim \omega_{\tilde{z}}'$.
Hence 
$\omega_z'\sim \omega_y$ 
is equivalent to $\omega_{\tilde{z}}'\sim \omega_{y}$.
From Theorem \ref{Thm:reductionequiv}(i),
this is equivalent that $|y_n|<1$ and $\tilde{z}=\tilde{y}$.
By definitions of $\tilde{z}$, $\tilde{y}$ and $\hat{y}$,
we can verify that this is equivalent that
$|y_n|<1$ and $z=\hat{y}$.
\qedh

\noindent
{\it Proof of Theorem \ref{Thm:finite}.}
From 
Lemma \ref{lem:reductionthree}(i), the statement holds.
\qedh
%
%
\sftt{Examples}
\label{section:fourth}
%
%
%
\ssft{Geometric progression states on $\co{2}$}
\label{subsection:fourthone}
In this subsection,
we show definitions and theorems 
 for \gprs s on $\co{2}$ of order $2$ and $\infty$
as examples of main theorems.
Let $s_1,s_2$ denote Cuntz generators of $\co{2}$.
%
%
\sssft{Case of order $2$}
\label{subsubsection:fourthoneone}
We summarize properties of geometric progression states on $\co{2}$ 
of order $2$.
Define $({\Bbb C}^3)_1
:=\{(w_1,w_2,w_3)\in {\Bbb C}^3:|w_1|^2+|w_2|^2+|w_3|^2=1\}$.
By definition,
$\omega$ is a geometric progression state on $\co{2}$ 
by $z=(z_1,z_2,z_3)\in ({\Bbb C}^3)_1$
if and only if $\omega$ satisfies the following equations:
%
%
\begin{equation}
\label{eqn:omesone}
\omega(s_1)=\overline{z_1},\quad 
\omega(s_2s_1)=\overline{z_2},\quad 
\omega(s_2^2)=\overline{z_3}.
\end{equation}
For a state $\omega$ on $\co{2}$ with 
the GNS representation $({\cal H},\pi,\omega)$,
the following are equivalent from 
Corollary \ref{cor:homoone}(i):
\begin{enumerate}
\item
$\omega$ is a geometric progression state by $z$.
\item
$\{z_1\pi(s_1)+z_2\pi(s_2s_1)+z_3\pi(s_2^2)\}\Omega=\Omega$.
\item
$\pi(s_1)^*\Omega=z_1\Omega$,
$\pi(s_2s_1)^*\Omega=z_2\Omega$
and 
$\pi(s_2^2)^*\Omega=z_3\Omega$.
\end{enumerate}
%
%
\begin{cor}
\label{cor:onfinite}
\begin{enumerate}
\item
A geometric progression state on $\co{2}$ by $z=(z_1,z_2,z_3)\in ({\Bbb C}^3)_1$
is unique if and only if $|z_3|<1$.
In this case, it is pure and 
we write such a state as $\omega_z'$.
\item
For $z,y\in {\cal W}_3:=\{(w_1,w_2,w_3)\in ({\Bbb C}^3)_1:|w_3|<1\}$,
$\omega_z'\sim \omega_y'$ if and only if $z=y$.
\item
For $z\in {\cal W}_3$,
$\omega_z'$ is equivalent to the Cuntz state $\omega_y$ on $\co{2}$
by $y=(y_1,y_2)\in ({\Bbb C}^2)_1$
if and only if $|y_2|<1$ and 
$z=(y_1,y_2y_1,y_2^2)$.
\end{enumerate}
\end{cor}
%
%
\pr
(i)
See Theorem \ref{Thm:mainb}(i).
(ii)
See Theorem \ref{Thm:mainthree}(i).
(iii)
See Theorem \ref{Thm:mequivalence}(iv).
\qedh

An idea of this study was brought from the following example.
Abe \cite{Abe} constructed a representation
$({\cal H},\pi)$
of $\co{2}$ with
a cyclic vector
$\Omega\in {\cal H}$
which satisfies
%
%
\begin{equation}
\label{eqn:abeexample}
\frac{1}{\sqrt{2}}\pi(s_{1}+s_{2}s_{1})\Omega=\Omega.
\end{equation}
By generalizing this,
we obtained a class of representations
of $\con$ and showed their properties (\cite{SE01}, Lemma 2.1).
From Corollary \ref{cor:onfinite},
the representation in (\ref{eqn:abeexample})
is unitarily equivalent to the GNS representation
by the geometric progression state 
$\omega_z'$ on $\co{2}$ for $z=(1/\sqrt{2},\,1/\sqrt{2},\,0)\in {\Bbb C}^3$
and it is irreducible.

%
%
\sssft{Case of order $\infty$}
\label{subsubsection:fourthonetwo}
A state $\omega$ on $\co{2}$ is a \gprs\ by $z=(z_1,z_2,\ldots)\in\ltno
:=\{(z_i):z_i\in {\Bbb C}\mbox{ for all }i\mbox{ and }
\sum_{i\geq 1}|z_i|^2=1\}$
if and only if 
$\omega$ satisfies
$\omega(s_2^{n-1}s_1)=\overline{z_n}$
for all $n\geq 1$.
From Corollary \ref{cor:homoone}(ii),
this is equivalent that 
$(z_1s_1+z_2s_2s_1+z_3s_2^2s_1+\cdots)\Omega=\Omega$,
or 
$\pi(s_2^{n-1}s_1)^*\Omega=z_n\Omega$ for all $n\geq 1$
where $({\cal H},\pi,\Omega)$ denotes
the GNS representation by $\omega$.
%
%
\begin{cor}
\label{cor:onitwo}
\begin{enumerate}
\item
For any $z\in\ltno$,
a \gprs\ on $\co{2}$ by $z$ exists uniquely and is pure.
We write it as $\omega_z'$.
\item
For $z,y\in\ltno$,
$\omega_z'\sim\omega_y'$
if and only if $z=y$.
\item
For $y=(y_1,y_2)\in({\Bbb C}^2)_1$,
let $\omega_y$ denote the Cuntz state on $\co{2}$ by $y$.
For $z\in \ltno$,
$\omega_z'\sim \omega_y$ if and only if 
$|y_2|<1$ and 
$z=(y_1,y_2y_1,y_2^2y_1,y_2^3y_1,\ldots)$.
In this case, 
$\omega_z'=\omega_y$.
\item
For $z\in\ltno$ and $y=(y_1,y_2,y_3)\in {\cal W}_{3}$,
$\omega_z'\sim \omega_y'$
if and only if 
%
%
\begin{equation}
\label{eqn:tildebbb}
z=(y_1,y_2,y_3y_1,y_3y_2,y_3^2y_1,y_3^2y_2,\ldots).
\end{equation}
In this case, 
$\omega_z'=\omega_y'$.
\end{enumerate}
\end{cor}
%
%
\pr
(i)
See Theorem \ref{Thm:mainb}(ii).
(ii)
See Theorem \ref{Thm:mainthree}(ii).
(iii) 
See Theorem \ref{Thm:reductionequiv}(i).
(iv) 
See Theorem \ref{Thm:reductionequiv}(ii).
\qedh

\noindent
From Corollary \ref{cor:onitwo}(iii),
$\omega_z'\sim \omega_y$ if and only if 
the sequence $z$ is just a geometric progression of complex numbers
with initial value $y_1$ and common ratio $y_2$
such that $|y_1|^2+|y_2|^2=1$ and $|y_2|<1$.
In this sense, we can regard geometric progression states
as generalizations of geometric progressions of complex numbers.

%
%
\ssft{Typical examples}
\label{subsection:fourthtwo}
%
%
\begin{ex}
\label{ex:first}
{\rm
Let $U(1):=\{c\in {\Bbb C}:|c|=1\}$.
As examples of 
Corollary \ref{cor:onfinite},
we show three kinds of states on $\co{2}$.
For $c\in U(1)$,
three states $\rho_c,\rho'_c,\rho''_c$ on $\co{2}$
which satisfy 
%
%
\begin{equation}
\label{eqn:rhoone}
\rho_c(s_1)=c,\quad 
\rho_c'(s_2s_1)=c,\quad 
\rho_c''(s_1+s_2s_1)=\sqrt{2}c
\end{equation}
exist uniquely and are pure.
They are the Cuntz state, the sub-Cuntz state
and the geometric progression state on $\co{2}$
by $(\overline{c},0)$,
$(0,\overline{c})\otimes (1,0)$ and 
$(\overline{c}2^{-1/2}, \overline{c}2^{-1/2},0)$,
respectively.
Any two distinct states in the set $\{\rho_c,\rho'_c,\rho''_c:c\in U(1)\}$
are not equivalent.
They can be unified as the geometric progression state 
$\rho_{c_1,c_2}$ which satisfies
%
%
\begin{equation}
\label{eqn:rhotwo}
\rho_{c_1,c_2}(c_1s_1+c_2 s_2s_1)=1\quad( (c_1,c_2)\in ({\Bbb C}^2)_1).
\end{equation}
In fact, we see that 
$\rho_c=\rho_{\overline{c},0}$,
$\rho'_c=\rho_{0,\overline{c}}$
and 
$\rho''_c=\rho_{\overline{c}2^{-1/2},\overline{c}2^{-1/2}}$ for 
$c\in U(1)$.
Any two distinct 
states in the set $\{\rho_{c_1,c_2}:(c_1,c_2)\in ({\Bbb C}^2)_1\}$
are not equivalent.
From this,
general geometric progression states
can be regarded as interpolations between two (or many)
special sub-Cuntz states of different orders.
}
\end{ex}

Next,
we show an example of Corollary \ref{cor:onitwo}.
Let $\zeta(x)$ denote the Riemann zeta function for
a positive real number $x>1$, that is,
$\zeta(x):=\sum_{n=1}^{\infty}n^{-x}$.
%
%
\begin{prop}
\label{prop:thensuch}
Fix $x>1$.
Assume that a state $\omega$ on $\co{2}$ satisfies
%
%
\begin{equation}
\label{eqn:zeta}
\sum_{n=1}^{\infty}
\frac{\omega(s_2^{n-1}s_1)}{n^{x/2}}=\sqrt{\zeta(x)}.
\end{equation}
Then the following holds.
\begin{enumerate}
\item
For any $x$, $\omega$ exists uniquely and is pure.
We write $\omega$ as $\kappa_x$.
Then
$\kappa_x\sim\kappa_{x'}$ if and only if $x=x'$.
\item
Let  $\rho_{c_1,c_2}$ be as in Example \ref{ex:first}.
For any $x$ and $(c_1,c_2)\in ({\Bbb C}^2)_1$,
$\kappa_x\not\sim \rho_{c_1,c_2}$.
\end{enumerate}
\end{prop}
%
%
\pr
(i)
Define $z(x):=(z_n(x))\in\ltn$ by
$z_n(x):=\{\zeta(x)n^{x}\}^{-1/2}\in {\Bbb R}$ for $n\geq 1$.
Then we see $\|z(x)\|^2=1$.
Hence $z(x)\in \ltno$.
From this and the definition, $\omega$ is the geometric progression
state on $\co{2}$ by $z(x)$.
From Corollary \ref{cor:onitwo}(i),
the uniqueness and purity hold.
From Corollary \ref{cor:onitwo}(ii),
the equivalence condition holds
because $z(x)=z(x')$ if and only if $x=x'$.

\noindent
(ii)
By definition,
we see that 
$\rho_{c_1,c_2}$ is the \gprs\ on $\co{2}$ by
$(c_1,c_2,0)\in {\cal W}_3$.
From Corollary \ref{cor:onitwo}(iv),
the statement holds.
\qedh

%
%
%
\ssft{Transformations of geometric progression states}
\label{subsection:fourththree}
We show transformations of geometric progression states of order $2$
and their equivalence.
Let $f$ be as in (\ref{eqn:gencc}) for $m=2n-1$.
Let $s_1,\ldots,s_n$ and $t_1,\ldots,t_{2n-1}$ denote
Cuntz generators of $\con$ and $\co{2n-1}$, respectively.
Define $\alpha\in {\rm Aut}\con$ 
and $\beta\in {\rm Aut}\co{2n-1}$ by 
%
%
\begin{equation}
\label{eqn:alphaone}
\alpha(s_i):=s_{n-i+1}\quad(i=1,\ldots,n),\quad 
\beta(t_j):=t_{2n-j}\quad(j=1,\ldots,2n-1).
\end{equation}
Then $\alpha^{-1}=\alpha$ and $\beta^{-1}=\beta$.
Define $f'\in {\rm Hom}(\co{2n-1},\con)$ by
$f':=\alpha\circ f\circ \beta$.
Then we can verify
%
%
\begin{equation}
\label{eqn:gti}
(f'(t_1),\ldots,f'(t_{2n-1}))=(s_1^2,s_1s_2,s_1s_3\ldots,s_1s_n,s_2,s_3,\ldots,s_n).
\end{equation}
By definition,
$\omega$ is an $f'$-sub-Cuntz state on $\con$ 
by $z=(z_1,\ldots,z_{2n-1})\in ({\Bbb C}^{2n-1})_1$
if and only if $\omega$ is a state on $\con$ which satisfies
%
%
\begin{equation}
\label{eqn:omegass}
\left\{
\begin{array}{l}
\omega(s_1^2)=\overline{z_1},\,\omega(s_1s_2)=\overline{z_2},\ldots,
\omega(s_1s_n)=\overline{z_n},\\
\\
\omega(s_2)=\overline{z_{n+1}},\ldots,
\omega(s_n)=\overline{z_{2n-1}}.
\end{array}
\right.
\end{equation}
%
%
%
\begin{prop}
\label{prop:geny}
Let $z=(z_1,\ldots,z_{2n-1})\in({\Bbb C}^{2n-1})_1$.
\begin{enumerate}
\item
For an $f'$-sub-Cuntz state $\omega$ on $\con$ by $z$, the following holds:
\begin{enumerate}
\item
$\omega\circ \alpha$ is the geometric progression state by 
$(z_{2n-1},z_{2n-2},\ldots,z_1)$.
\item
$\omega$ is unique if and only if $|z_1|<1$.
In this case, $\omega$ is pure and we write $\omega$ as $\eta_z$.
\end{enumerate}
\item
Define 
${\cal W}'_{2n-1}:=\{(w_1,\ldots,w_{2n-1})\in 
({\Bbb C}^{2n-1})_1: |w_1|<1\}$.
For $z,y\in {\cal W}'_{2n-1}$,
$\eta_z\sim\eta_y$ if and only if $z=y$.
\end{enumerate}
\end{prop}
%
%
\pr
(i)
(a) By assumption,  the statement holds.

\noindent
(b)
From (a) and Theorem \ref{Thm:mainb}(i),
the statement holds.

\noindent
(ii)
From (i)(a) and Theorem \ref{Thm:mainthree}(i),
the statement holds.
\qedh


\appendix \section*{Appendix} 
%
%
\sftt{Proof of Lemma \ref{lem:words}}
\label{section:appone}
In order to prove Lemma \ref{lem:words},
we show a lemma associated with a free semigroup
and its subsemigroups.
Fix $2\leq n<\infty$.
Let $s_1,\ldots,s_n$ denote Cuntz generators of $\con$.
Let ${\cal I}_n$ be as in (\ref{eqn:idefinition})
and ${\sf S}_n:=\{s_J: J\in {\cal I}_n,\,J\ne\emptyset\}$.
Then ${\sf S}_n$ is 
the non-selfadjoint subsemigroup of $\con$ generated by $s_1,\ldots,s_n$ 
without unit, and it is free \cite{Howie1976,SO}.
For two nonempty subsets $X,Y\subset {\sf S}_n$,
define $XY:=\{xy: x\in X,\,y\in Y\}$
and $X^0:=\{I\}$, $X^a:=X^{a-1}X$ for $a\geq 1$.
Let $A$ and $B$ denote subsemigroups of ${\sf S}_n$
generated by $\{s_1,\ldots,s_{n-1}\}$
and $\{s_n\}$, respectively:
%
%
\begin{equation}
\label{eqn:abset}
A:=\langle \{s_1,\ldots,s_{n-1}\}\rangle,\quad 
B:=\langle\{s_n\}\rangle=\{s_n^a:a\geq1 \}.
\end{equation}
Since ${\sf S}_n$ coincides with the free product $A*B$ \cite{Howie1976,SO},
any $x\in {\sf S}_n$ belongs to one of the following:
%
%
\begin{equation}
\label{eqn:abab}
A(BA)^a,\quad 
A(BA)^aB,\quad 
(BA)^{a'},\quad 
(BA)^aB
\end{equation}
for some $a\geq 0$ and $a'\geq 1$.

By identifying $t_i$ with $f(t_i)$ in Definition \ref{defi:lexico},
we regard $t_1,t_2,\ldots\in \con$.
Let ${\sf T}_m$ denote the subsemigroup of ${\sf S}_n$
generated by $t_1,t_2,\ldots$.
Then ${\sf T}_m$ is also free.
Since $t_i=s_i$ for $i=1,\ldots,n-1$, $A\subset {\sf T}_m$.
%
%
\begin{lem}
\label{lem:wordone}
For any $2\leq m\leq \infty$ and $a\geq 1$, $(BA)^a\subset {\sf T}_m$.
\end{lem}
%
%
\pr
It is sufficient to show $BA\subset {\sf T}_m$.

\noindent
(i) Assume $m<\infty$.
If $x\in BA$, then $x=s_n^a y$ for some $a\geq 1$ and $y\in A$.
Then we can write $a=a'k+a''$ and 
$y=s_iy'$ for some $a'\geq 0$, $0\leq a''\leq k-1$, $(a',a'')\ne (0,0)$,
and some $1\leq i\leq n-1$ and $y'\in A\cup \{I\}$.
From these and $A\subset {\sf T}_m$,
$x=s_n^{a'k+a''}(s_iy')=t_m^{a'}t_{(n-1)a''+i}y' \in {\sf T}_m$.

\noindent
(ii) Assume $m=\infty$.
If $x\in BA$, then $x=s_n^a y$ for some $a\geq 1$ and $y\in A$.
Then we can write $y=s_iy'$ for some $1\leq i\leq n-1$ and $y'\in A\cup \{I\}$.
From these and $A\subset {\sf T}_{\infty}$,
$x=s_n^{a}s_iy'=t_{(n-1)a+i}y' \in {\sf T}_{\infty}$.
\qedh

\noindent
{\it Proof of Lemma \ref{lem:words}.}
(i) Assume $m<\infty$.
Let ${\sf T}:={\sf T}_m$ and $x=s_J$. 
If $x\in {\sf T}$,
then the statement holds because ${\sf T}={\sf T}s_n^0$.
We check the statement for each case in (\ref{eqn:abab}) as follows.
\begin{enumerate}
\item[(a)]
Assume $x\in A(BA)^a$ for some $a\geq 0$.
Since $A\subset {\sf T}$,
$x\in A(BA)^a\subset {\sf T}$
from Lemma \ref{lem:wordone}.
Hence the statement holds.
\item[(b)]
Assume $x\in A(BA)^aB$ for some $a\geq 0$.
From $A\subset {\sf T}$ and Lemma \ref{lem:wordone},
we can write $x=x's_n^b$ for some $x'\in {\sf T}$ and $b\geq 1$.
Then we can write $b=b'k+b''$ for some
$b'\geq 0$, $0\leq b''\leq k-1$ and  $(b',b'')\ne (0,0)$.
Then $x=x's_n^{b'k+b''}=x't_m^{b'}s_n^{b''}\in {\sf T}s_n^{b''}$.
Hence the statement holds.
\item[(c)]
Assume $x\in (BA)^{a'}$ for some $a'\geq 1$.
From Lemma \ref{lem:wordone}, the statement holds.
\item[(d)]
Assume $x\in (BA)^aB$ for some $a\geq 0$.
Along with (b),  the statement holds.
\end{enumerate}

\noindent
(ii)
Let ${\sf T}:={\sf T}_{\infty}$.
Remark that 
if $x\in {\sf T}$,
then the statement holds because ${\sf T}={\sf T}s_n^0$.
Except the case of (b) in the proof of (i),
the statement for each case in (\ref{eqn:abab}) holds as (i).
Assume $x\in A(BA)^aB$ for some $a\geq 1$.
From $A\subset {\sf T}$ and Lemma \ref{lem:wordone},
we can write $x=x's_n^b$ for some $x'\in {\sf T}$
and $b\geq 1$.
Hence the statement holds.
\qedh

%
%

%
\label{Lastpage}
\end{document}